\documentclass{article}
\usepackage{mathrsfs}
\usepackage{ifthen}
\usepackage{amssymb}
\usepackage[all]{xy}
\usepackage[leftbars]{changebar}
\usepackage{vmargin}

\setpapersize{A4}
\setmarginsrb{37mm}{20mm}{30mm}{37mm}{10mm}{6mm}{0mm}{10mm}

\def\institute#1{\gdef\@institute{#1}}

\newcommand{\tiret}{\rule[0.6ex]{1.3ex}{0.26ex}}
\newcommand{\defmot}[1]{\emph{#1}}
\newcommand{\ins}[1]{{\ensuremath{\textrm{#1}}}}
\newcommand{\qetq}{\quad\textrm{and}\quad}
\newcommand{\qqetqq}{\qquad\textrm{and}\qquad}
\newcommand{\donne}[2]{\noindent\ensuremath{#1\Rightarrow #2.}\quad}
\newcommand{\pa}[1]{\ensuremath{\left(#1\right)}}
\newcommand{\trans}[1]{\ensuremath{\vphantom{A}{^t}{}%
        {\mathop{\mathstrut #1 \!}}}\,}
\newcommand{\wt}{\widetilde}
\newcommand{\wh}{\widehat}
\newcommand{\im}{\ensuremath{\mathrm{Im\,}}}
\newcommand{\tr}{\ensuremath{\mathrm{tr\,}}}
\newcommand{\accolades}[1]{\ensuremath{\left\{#1\right\}}}
\newcommand{\paa}[1]{\ensuremath{\accolades{#1}}}

\newcommand{\La}{\Lambda}
\newcommand{\Om}{\Omega}
\newcommand{\De}{\Delta}
\newcommand{\al}{\alpha}
\newcommand{\dz}{\zeta}
\newcommand{\om}{\omega}
\newcommand{\ep}{\varepsilon}
\newcommand{\de}{\delta}
\newcommand{\la}{\lambda}
\newcommand{\ga}{\gamma}
\newcommand{\ph}{\varphi}
\newcommand{\si}{\sigma}
\newcommand{\be}{\beta}

\newcommand{\Sscr}{\mathscr{S}}
\newcommand{\Bscr}{\mathscr{B}}
\newcommand{\Gscr}{\mathscr{G}}
\newcommand{\Fscr}{\mathscr{F}}
\newcommand{\Nscr}{\mathscr{N}}
\newcommand{\Hyp}{\mathscr{H}}

\newcommand{\NN}{\mathbb{N}}
\newcommand{\ZZ}{\mathbb{Z}}
\newcommand{\UU}{\mathbb{U}}
\newcommand{\QQ}{\mathbb{Q}}
\newcommand{\CC}{\mathbb{C}}
\newcommand{\FF}{\mathbb{F}}

\newcommand{\Sfr}{\mathfrak{S}}

\newcommand{\GL}[1]{\ensuremath{\mathrm{GL}}(#1)}
\newcommand{\id}{\ensuremath{\mathrm{id}}}
\newcommand{\ra}{\rightarrow}
\newcommand{\Res}[3]{\ensuremath{\mathrm{Res}_{#1}^{#2}(#3)}}
\newcommand{\intnn}[2]{\ensuremath{[\![ \, #1 \,,\, #2 \,]\!]}}
\newcommand{\homm}[4]{\ensuremath{\mathrm{Hom}_{#1}^{#2}(#3,#4)}}
\newcommand{\Hom}[3]{\homm{#1}{}{#2}{#3}}
\newcommand{\SSI}{\ensuremath{\;\Longleftrightarrow\;}}
\newcommand{\mat}[2]{\ensuremath{\mathrm M}_{#1}(#2)}
\newcommand{\gal}[2]{\ensuremath{\mathrm{Gal}}(#1/#2)}
\newcommand{\galq}{\gal{\overline{\QQ}}{\QQ}}
\newcommand{\Ende}[3]{\ensuremath{\mathrm{End}_{#1}^{#2}(#3)}}
\newcommand{\End}[2]{\Ende{#1}{}{#2}}
\newcommand{\diag}{\ensuremath{\mathrm{diag\,}}}

\newcommand{\Bsdchi}{\ensuremath{B_\si(d,\chi)}}
\newcommand{\bsdchigamma}{\ensuremath{b_{\si,\gamma}(d,\chi)}}
\newcommand{\Bsdchigamma}{\ensuremath{B_{\si,\gamma}(d,\chi)}}
\newcommand{\bsdchigammastar}{\ensuremath{b^{*}_{\si,\gamma}(d,\chi)}}
\newcommand{\Bsdchigammastar}{\ensuremath{B^*_{\si,\gamma}(d,\chi)}}

\newcommand{\fonction}[5]{%
        \ensuremath{#1\colon
        \left\{\hskip -1.5 mm
        \begin{array}{c@{\ }c@{\ }l}
        \medskip #2 & \longrightarrow & #3 \\
        #4 & \longmapsto & #5 \\
        \end{array}
        \right .
        }}

\newcommand{\Somd}[2]{\ensuremath{\bigoplus\limits}_{#1}^{#2}}
\newcommand{\Somdd}[1]{\ensuremath{\bigoplus\limits}_{#1}^{}}
\newcommand{\Sum}[2]{\ensuremath{\textstyle{\sum\limits_{#1}^{#2}}}}
\newcommand{\Prod}[2]{\ensuremath{\prod\limits}_{#1}^{#2}}
\newcommand{\Prodd}[1]{\ensuremath{\prod\limits}_{#1}}
\newcommand{\Sumu}[1]{\Sum{#1}{}}

\newenvironment{rs}{\begin{list}{\tiret~}{
 \topsep=0.3ex
 \itemsep=0.3ex
 \labelsep=0em
 \parsep=0em
 \listparindent=1em
 \itemindent=0em
 \settowidth{\labelwidth}{--~}
 \leftmargin=\labelwidth
}}{\end{list}}

\newcommand{\numtoi}[1]{
  \ifthenelse{ \equal{#1}{1} }{i}{
  \ifthenelse{ \equal{#1}{2} }{ii}{
  \ifthenelse{ \equal{#1}{3} }{iii}{
  \ifthenelse{ \equal{#1}{4} }{iv}{
  \ifthenelse{ \equal{#1}{5} }{v}{
  \ifthenelse{ \equal{#1}{6} }{vi}{
  \ifthenelse{ \equal{#1}{7} }{vii}{
  \ifthenelse{ \equal{#1}{8} }{viii}{
  \ifthenelse{ \equal{#1}{9} }{ix}{
  \ifthenelse{ \equal{#1}{10} }{x}{
  \ifthenelse{ \equal{#1}{11} }{xi}{
  \ifthenelse{ \equal{#1}{12} }{xii}{
  \ifthenelse{ \equal{#1}{13} }{xiii}{
  \ifthenelse{ \equal{#1}{14} }{xiv}{
  \ifthenelse{ \equal{#1}{15} }{xv}{
  \ifthenelse{ \equal{#1}{16} }{xvi}{
  \ifthenelse{ \equal{#1}{17} }{xvii}{
  \ifthenelse{ \equal{#1}{18} }{xviii}{
    ERREUR,~MODIFIER~LA~MACRO~numtoi
  } } } } } } } } } } } } } } } } } }
}

\newcounter{cretraiti}
\newenvironment{ri}[1]{
\begin{list}{($\numtoi{\thecretraiti}$)}{
  \usecounter{cretraiti}
  \topsep=0.5ex
  \itemsep=0.3ex
  \labelsep=0.3em
  \parsep=0ex
  \listparindent=1em
  \settowidth{\labelwidth}{(#1)}
  \leftmargin=\labelwidth
}}{\end{list}
}

\newenvironment{rn}[1]{\begin{list}{}{
 \topsep=0.5ex
 \itemsep=0.3ex
 \labelsep=0.3em
 \parsep=0ex
 \listparindent=1em
 \settowidth{\labelwidth}{(#1)}
 \leftmargin=\labelwidth
}}{\end{list}}

\newenvironment{rcl}[1][r]%
    {\ensuremath{\begin{array}[t]{#1@{\ }c@{\ }l}}}%
    {\end{array}}
\newenvironment{accolade}[1][c]%
    {\ensuremath{\left \{ \hskip -1.5 mm \begin{array}{#1@{\quad}l}}}%
    {\end{array}\right.}
\newenvironment{matrice}%
{\ensuremath{\left[\begin{array}{c@{\ }c}}}%
{\end{array}\right]}

\newcommand{\spacebeforeenv}{\vspace{1ex}}
\newcommand{\spaceafterenv}{\vspace{1ex}}

\newenvironment{traitsurlecote}
{\cbstart
\setcounter{changebargrey}{0}      
}
{\cbend
}

\newcommand{\mysquare}{\rule[0.3ex]{1ex}{1ex}}
\newcommand{\myqed}{\hfill \mysquare}
\newcommand{\myqedenv}{}
\newcommand{\smartqed}{\myqed}

\newcommand{\envfont}{\sf\bfseries}

\newcounter{ctheo}
\renewcommand{\thectheo}{\arabic{ctheo}}

\newcommand{\metaenvironnementhm}[2]  
{
  \refstepcounter{ctheo}
  \ifthenelse{ \equal{#2}{toto} }{
    \spacebeforeenv \begin{traitsurlecote} \noindent {\envfont #1 \thectheo}
  }{
    \spacebeforeenv \begin{traitsurlecote} \noindent {\envfont #1 \thectheo{} \tiret\, #2.}
  }
}

\newenvironment{theorem}[1][toto]
{
  \metaenvironnementhm{Theorem}{#1}
}
{\end{traitsurlecote}\spaceafterenv}

\newenvironment{proposition}[1][toto]
{
  \metaenvironnementhm{Proposition}{#1}
}
{\end{traitsurlecote}\spaceafterenv}

\newenvironment{lemma}[1][toto]
{ \metaenvironnementhm{Lemma}{#1} }
{\end{traitsurlecote}\spaceafterenv}

\newenvironment{corollary}[1][toto]
{ \metaenvironnementhm{Corollary}{#1} }
{\end{traitsurlecote}\spaceafterenv}

\newenvironment{definition}[1][toto]
{ \metaenvironnementhm{Definition}{#1} }
{\end{traitsurlecote}\spaceafterenv}

\newenvironment{notation}[1][toto]
{ \metaenvironnementhm{Notation}{#1} }
{\end{traitsurlecote}\spaceafterenv}

\newenvironment{hypothesis}[1][toto]
{ \metaenvironnementhm{Hypothesis}{#1} }
{\end{traitsurlecote}\spaceafterenv}

\newenvironment{proofproof}
{\noindent  {\envfont Proof.~}}
{\spaceafterenv}

\newenvironment{keyword}
{\noindent  {\envfont Keywords.~}}
{\spaceafterenv}

\newenvironment{msc}
{\noindent  {\envfont Mathematics Subject Classification (2000).~}}
{\spaceafterenv}

\newcommand{\metaenvironnement}[2]  
{
  \refstepcounter{ctheo}
  \ifthenelse{ \equal{#2}{toto} }{
    \spacebeforeenv \noindent {\envfont #1 \thectheo}
  }{
    \spacebeforeenv \noindent {\envfont  #1 \thectheo{} \tiret\, #2.}
  }
}

\newcommand{\metaenvironnementsimple}[2]  
{
  \ifthenelse{ \equal{#2}{toto} }{
    \spacebeforeenv \noindent {\envfont #1}      
  }{
    \spacebeforeenv \noindent {\envfont #1 \tiret\, #2.}
  }
}

\newenvironment{remark}[1][toto]
{ \metaenvironnement{Remark}{#1} }
{\myqedenv\spaceafterenv}

\begin{document}

\title{Exterior Algebra Structure for Relative Invariants of Reflections Groups}

\author{Beck Vincent}

\maketitle

\begin{abstract}
Let $G$ be a reflection group acting on a vector space $V$ (over a field with zero characteristic).
We denote by $S(V^*)$ the coordinate ring of $V$, by $M$ a finite dimensional $G$-module and by $\chi$ a
one-dimensional character of $G$.
In this article, we define an algebra structure on the isotypic component associated to $\chi$ of the algebra
$S(V^*) \otimes \Lambda(M^*)$. This structure is then used to obtain various generalizations of usual criterions on regularity of integers.
\end{abstract}


\begin{msc} 13A50, 15A75  \end{msc}


\section{Introduction}\label{intro}

In the first part of this article, we will study the following situation. Let $G$ be a reflection group acting on the vector space~$V$,
$M$ be a finite dimensional representation of $G$ and $\chi$ be a one-dimensional character of $G$. Following the ideas of Shepler~\cite{shepler},
we construct an exterior algebra structure on the $\chi$-isotypic component of $T^{-1}S(V^*) \otimes \La(M^*)$ for a suitable multiplicative set
$T$ of $S(V^*)$. This work is in line with the articles~\cite{orlik-solomon}, \cite{shepler}, \cite{lehrer} and~\cite{beck} which construct
algebra structures on the $\chi$-isotypic of the algebra $S(V^*) \otimes \La(M^*)$ under conditions over the restrictions of $M$ and $\chi$ to certain subgroups of $G$.
Here, the idea is to transfer the hypotheses on $M$ and $\chi$ to conditions on the base ring :
we substitute $S(V^*)$ in a bigger ring (a fraction ring of $S(V^*)$) in which some linear forms associated to hyperplanes of $G$ are invertible.
The conditions will be held by the "bad" hyperplanes that are needed to be invertible.
The main results of~\cite{orlik-solomon}, \cite{shepler}, \cite{lehrer} et~\cite{beck} are exceptional cases
of proposition~\ref{prop-conclusion-alg-exterieure} (see remark~\ref{rem-shepler-os}).
The article~\cite{shepler-hartmann} explains the situation in prime characteristic.

In the second part of this article, we will give consequences of the exterior algebra structure with links to the notion of regular integers.
These consequences are similar to those that can be found in~\cite{orlik-solomon}, \cite{lehrer-michel}, \cite{lehrer} and~\cite{bonnafe-lehrer-michel}.

Various types of hyperplanes appear in the first part of the article. The hyperplanes to invert (the multiplicative set $T$) are chosen following these types.
The third part studies these types for concrete reflection groups : the symmetric group, $G(de,e,2)$, $G(d,1,r)$ and the exceptional group $G_4$, $G_5$ and $G_{24}$.

Let us begin with some usual definitions and notations.

\begin{definition}[Reflection] Let $k$ be a field of characteristic $0$ and $V$  a finite dimensional vector space over $k$.
Any $g \in \GL{V}$ so that $g$ is of finite order and $\ker(g - 1)$ is an hyperplane of $V$ is called a \defmot{reflection}.
\end{definition}

\begin{definition}[Reflection Groups]\label{dfn-grc} Let $k$ be a field of characteristic $0$.
$(G,V)$ is said to be a \defmot{reflection group over $k$} if $V$ is a finite dimensional vector space over $k$ (we denote by $\ell$ the dimension of $V$)
and $G$ be a finite subgroup of $\GL{V}$ generated by reflections.
It will often be more comfortable to write "let $G$ be a reflection group" omitting the vector space $V$.
\end{definition}

\begin{notation}[Reflection and hyperplane]\label{notation-hyperplan-reflexion}
Let $(G,V)$ be a reflection group. We denote by $\Sscr$ the set of reflections of $G$ and $\Hyp$ the set of hyperplanes of $G$ :
$$\Sscr = \{s \in G,\quad \dim \ker(s- \id)= \dim V -1 \} \qquad \ins{and } \qquad \Hyp = \{ \ker(s -\id),\quad s \in \Sscr \}\,.$$
\end{notation}

\begin{notation}[Around a hyperplane]\label{notation-autour-hyperplan} Let $(G,V)$ be a reflection group. For $H \in \Hyp$,
\begin{rs}
\item one chooses $\al_H \in V^*$ a linear form with kernel $H$;

\item one sets $G_H = \ins{Fix}_{G}(H) = \{g \in G, \quad \forall\,x \in H, \quad gx = x \}$.
      This is a cyclic subgroup of $G$. We denote by $e_H$ its order and by $s_H$ its generator with determinant $\zeta_H = \exp(2i\pi/e_H)$;

\item For any finite dimensional $kG$-module $N$ and for $j \in \intnn{0}{e_H-1}$,
      we define the integers $n_H(N)$ and $n_{j,H}(N)$ by
      $$\Res{G_H}{G}{N} = \Somd{j=0}{e_H-1} n_{j,H}(N) \det{}^{-j}\qquad \ins{and} \qquad
          n_H(N) = \Sum{j=0}{e_H-1}j n_{j,H}(N)\,;$$
      \noindent the integer $n_{j,H}(N)$ is nothing else but the multiplicity of ${\zeta_H}^j$ as an eigenvalue of $s_H$ acting on $N^*$;

\item let $\chi : G \ra k^\times$ be a linear character of $G$, we denote by $k_\chi$ the representation of $G$ with character $\chi$
      over $k$ and $n_H(\chi)$ for $n_H(k_\chi)$; by definition, $n_H(\chi)$ is the unique integer $j$ verifying
      $0 \leq j < e_H$ and $\chi(s_H) = \det(s_H)^{-j}$\,. Finally, for any $kG$-module $N$, we denote by $N^\chi=\{x \in N, \quad gx=\chi(g)x\}$
      the $\chi$-isotypic component of $N$.
\end{rs}
\end{notation}

\vskip1ex
\begin{definition}[Polynomial function associated to a representation]\label{dfn-polynome-representation}
Let $(G,V)$ be a reflection group and $N$ be a finite dimensional $G$-module. We set
$$Q_N = \Prodd{H \in \Hyp} {\al_H}^{n_H(N)} \in S(V^*)\,.$$
\noindent When $\chi$ is a linear character of $G$, we set $Q_\chi$ rather than $Q_{k_\chi}$, so that
$$Q_{\chi} = \Prodd{H \in \Hyp}{\al_H}^{n_H(\chi)} \in S(V^*).$$
\end{definition}

\section{Construction of the Algebra Structure}\label{sec-construction-alg-exterieure}

In this section and the next one, we fix a reflection group $(G,V)$ over $k$, a $kG$-module $M$ with dimension $r$ and $\chi: G \ra k^\times$ a linear character of~$G$.
We denote by $\det_M$ (resp. $\det_{M^{*}}$) the determinant of the representation $M$ (resp. $M^*$).

\begin{notation}
Let $\Bscr \subset \Hyp$ be a $G$-stable subset (which will be the "bad" hyperplanes), we denote by $\Gscr= \Hyp \setminus \Bscr$
and
$T = \langle \al_H,\quad H \in \Bscr \rangle_{\ins{mult. set}}$ the multiplicative subset of $S(V^*)$ associated to $\Bscr$.
We then set $\Om = T^{-1}S(V^*) \otimes \La(M^*)$ and $\Om^p = T^{-1}S(V^*) \otimes \La^p(M^*)$ for $p \in \intnn{0}{r}$.
Thus we have
$$\Om^\chi = \Somd{p=0}{r}\, (\Om^p)^\chi\,.$$
\end{notation}

As said in the introduction, the idea of the roof is to bring the "bad" hyperplanes together in the subset $\Bscr$.
Thus in subsection~\ref{ssec-invrel-alg-ext-hyperplan}, we begin to define some "hyperplane types" so that we are able to differentiate the hyperplanes of $\Hyp$.
In subsection~\ref{ssec-invrel-alg-ext-struct-alg}, we construct the algebra structure on $\Om^\chi$.
Finally, in subsection~\ref{ssec-invrel-algebre-exterieure}, we study this algebra structure.

\subsection{Hyperplanes}\label{ssec-invrel-alg-ext-hyperplan}

To define the notion of $(M,\chi)$-acceptable hyperplane, we need to introduce a notation which will also be useful for the next subsection.

\begin{notation}\label{notation-ji} For $H \in \Hyp$, we denote by $(j_1, \ldots, j_r)$ a family of integers such that for all $i \in \intnn{1}{r}$
we have $0 \leq j_i \leq e_H-1$ and $\ins{Res}_{G_H}^{G}(M) = \det^{-j_1} \oplus \cdots \oplus \det^{-j_r}$. The family $(j_1, \ldots, j_r)$ is not unique
but so is the multi-set associated to it. The integers $j_i$ are closely related to the integers $n_{j,H}(M)$ (see notation~\ref{notation-autour-hyperplan}).
Precisely, for $j \in \intnn{0}{e_H-1}$, $n_{j,H}(M)$ is the number of $i \in \intnn{1}{r}$ so that $j_i = j$, so that we have
$$\Sum{i=1}{r} j_i = \Sum{j=0}{e_H -1} j n_{j,H}(M) = n_H(M)\,.$$
\noindent In addition, the family $(j_1,\ldots, j_r)$ is so that the eigenvalues of $s_H$ acting on $M^*$ are the ${\dz_H}^{j_i}$ for $i \in \intnn{1}{r}$.
\end{notation}

We define four types of hyperplanes (two types associated with a linear representation of $G$ and two others associated with a couple constituted of a representation of
$G$ and a one-dimensional character of $G$). The $M$-excellent hyperplanes or the $(M,\chi)$-good hyperplanes will be those that we do not need to inverse
(see hypotheses~\ref{hyp-excellent-hyperplan} and~\ref{hyp-action-reflexion}).

\begin{definition}[Hyperplane types]\label{dfn-type-hyperplan} Let $H \in \Hyp$; $H$ is said to be
\begin{ri}{(iii)}
\item \defmot{$M$-good} if $n_H(M) < e_H$ and \defmot{$M$-bad} else;
\item \defmot{$M$-excellent} if $s_H$ acts on $M$ as a reflection;
\item \defmot{$(M,\chi)$-good} if $n_H(M) + n_H(\chi) < e_H$;
\item \defmot{$(M,\chi)$-acceptable} if for all partition of the set $\intnn{1}{r}$ in two disjoint sets (denoted respectively by $I_1$ and $I_2$), we have
$$e_H - n_H(\chi) > \Sumu{i \in I_1} j_i \qquad\ins{or}\qquad e_H - n_H(\chi) > \Sumu{i \in I_2} j_i\,.$$
\end{ri}
\end{definition}

In the next remark, we study the links between the preceding notions of hyperplane types.

\begin{remark}[Good and excellent hyperplanes]\label{rem-bon-excellent-hyperplan}
Let $H \in \Hyp$. We denote by ${1}$ the trivial character on $G$. Let show the following properties.

\begin{ri}{(iii)}
\item $H$ is $M$-excellent if and only if $H$ is $M$-good and $M^*$-good.
\item If $H$ is $(M,\chi)$-good then $H$ is $M$-good.
\item $H$ is $M$-excellent if and only if $H$ is $(M,\chi)$-acceptable for all $\chi \in \Hom{\ins{gr}}{G}{k^{\times}}$.
\item If $H$ is $(M,\chi)$-good then $H$ is $(M,\chi)$-acceptable.
\item $H$ is $M$-good if and only if $H$ is $(M,1)$-good.
\item If $H$ is $M$-good then $H$ is $(M,1)$-acceptable.
\end{ri}

Let us begin with~$(i)$. We start to express $n_H(M^*)$ using $n_H(M)$ :
$$n_H(M^*) = \Sum{j=1}{e_H-1} (e_H-j)n_{j,H}(M) = e_H (r - n_{0,H}(M)) - n_H(M)\,.$$
\noindent We have $n_H(M) = 0$ if and only if $n_H(M^*) = 0$ if and only if $s_H$ acts trivially on $M$. We then deduce
$$n_H(M) \leq e_H -1 \quad \ins{and} \quad n_H(M^*) \leq e_H -1 \qquad \SSI \qquad r - n_{0,H}(M) <2\,.$$
\noindent Since $n_{0,H}(M)$ (resp. $n_{0,H}(M^*)$) is the multiplicity of $1$ as an eigenvalue of $s_H$ acting on $M$
(resp. $M^*$), we obtain $n_{0,H}(M) = n_{0,H}(M^*)$ and the condition $n_{0,H}(M) \in \{r-1,r\}$ can be expressed geometrically
as $s_H$ acts trivially on $M$ or acts as a reflection on $M$.
In particular, an $M$-excellent hyperplane is always $M$-good.

Let us show~$(ii)$. We have $n_H(\chi)\geq 0$, so an $(M,\chi)$-good hyperplane is always $M$-good.

Now, let us consider~$(iii)$. Let us assume that $H$ is an $M$-excellent hyperplane. There exists at most one $i_0 \in \intnn{1}{r}$ so that $j_{i_0}$ is nonzero.
When $I_1$ and $I_2$ are two disjoint subsets of $\intnn{1}{r}$, only one of those two sets can contain $i_0$.
Thus, we have
$$\Sumu{i \in I_1} j_i = 0 < e_H - n_H(\chi) \qquad \ins{or} \qquad \Sumu{i \in I_2} j_i = 0 < e_H - n_H(\chi)\,.$$
\noindent Reciprocally, let us assume that $H$ is not $M$-excellent. We then deduce that there exists $i_1 \neq i_2$ so that $j_{i_1} \neq 0$ and $j_{i_2} \neq 0$.
In addition, we know the existence of a linear character $\chi$ of $G$ so that $n_H(\chi)= e_H-1$ by Stanley's theorem~\cite{stanley}.
The disjoint sets $I_1=\{i_1\}$ and $I_2=\{i_2\}$ verify
$$\Sumu{i \in I_1} j_i \geq 1 = e_H - n_H(\chi) \qquad \ins{and} \qquad \Sumu{i \in I_2} j_i \geq 1 = e_H - n_H(\chi)\,.$$
\noindent We then deduce that $H$ is not $(M,\chi)$-acceptable.

Let us show~$(iv)$. Let $I$ be a subset of $\intnn{1}{r}$ so that
$$e_H - n_H(\chi) \leq \Sumu{i \in I} j_i\,.$$
\noindent Any such $I$ contains every $i \in \intnn{1}{r}$ so that $j_i \neq 0$ since
$$\Sumu{i \in \intnn{1}{r}} j_i = n_H(M) \leq e_H - n_H(\chi)\,.$$
\noindent Thus, any set $I'$ disjoint of $I$ verifies
$$\Sumu{i \in I'} j_i = 0 < e_H - n_H(\chi)\,.$$

Let us show~$(v)$ and $(vi)$. Since $n_H(1) = 0$ for all $H \in \Hyp$, an hyperplane is $M$-good if and only if it is $(M,1)$-good.
$(iv)$ shows that such an hyperplane is $(M,1)$-acceptable. \smartqed
\end{remark}

Since lots of hyperplanes of reflections groups verify $e_H= 2$ (for example this is the case for Coxeter groups but not only), we focus on this specific case.

\begin{remark}[Hyperplanes with ${e_H=2}$]\label{rem-type-hyperplan-ordre2}
Let $H \in \Hyp$ so that $e_H=2$.

Then $H$ is $M$-good if and only if $H$ is $M$-excellent
(that is if the multiplicity of the eigenvalue $-1$ of $s_H$ acting on $M$ is not bigger than one).

If $\chi(G_H) \neq 1$ then
\begin{ri}{(ii)}
\item $H$ is $(M,\chi)$-acceptable if and only if $s_H$ acts on $M$ as a reflection or acts trivially on $M$
(that is if the multiplicity of the eigenvalue $-1$ of $s_H$ acting on $M$ is not bigger than $1$);
\item $H$ is $(M,\chi)$-good if and only if $s_H$ acts trivially on $M$ (that is if the multiplicity of the eigenvalue $-1$ of $s_H$ acting on $M$ is zero);
\end{ri}

If $\chi(G_H) = 1$ then
\begin{ri}{(ii)}
\item $H$ is $(M,\chi)$-acceptable if and only if the multiplicity of the eigenvalue $-1$ of $s_H$ acting on $M$ is not bigger than $3$;
\item $H$ is $(M,\chi)$-good if and only if $s_H$ acts on $M$ as a reflection or acts trivially on $M$
(that is if the multiplicity of the eigenvalue $-1$ of $s_H$ acting on $M$ is not bigger than $1$).
\end{ri}

Following the remark~\ref{rem-bon-excellent-hyperplan}, it is enough to show that if $H$ is $M$-good then $H$ is $M$-excellent.
By hypothesis, we have $n_H(M) =n_{1,H}(M) < 2$ and thus $n_{0,H}(M) = r-n_{1,H}(M) \in \{r,r-1\}$.

Let us assume that $\chi(G_H)\neq 1$. We have $n_{H}(\chi) \neq 0$ and then $n_H(\chi) = 1$. The definition of "acceptability" shows us that
if $H$ is $(M,\chi)$-acceptable then $H$ is $(M,\chi')$-acceptable for all linear characters $\chi'$ verifying $n_H(\chi') \leq 1$.
But every linear character $\chi'$ verifies $n_H(\chi') \leq 1$, thus we have $H$ is $(M,\chi')$-acceptable for all linear characters $\chi'$ of $G$.
The remark~\ref{rem-bon-excellent-hyperplan} shows that $H$ is $M$-excellent.
In addition, by definition of $(M,\chi)$-good, $H$ is $(M,\chi)$-good if and only if $n_H(M)<1$ that is $n_H(M) = 0$.

Let us assume that $\chi(G_H)= 1$. We have $n_H(\chi) = 0$ and then $H$ is $(M,\chi)$-good if and only if $n_H(M)<2$ (that is $H$ is $M$-good).
In addition, in our case, we have $j_i \in \{0,1\}$, the multiplicity of $-1$ as an eigenvalue of $s_H$ is the number of $i$ so that $j_i \neq 0$.
If they are more than $4$, we can divide them in two sets of two elements and the hyperplane is not $(M,\chi)$-good.
If they are not more than $4$, two disjoints set of $\intnn{1}{r}$ cannot both contain
two integers $i$ so that $j_i \neq 0$ and finally $H$ is $(M,\chi)$-acceptable. \smartqed
\end{remark}

\subsection{Construction of an Algebra Structure}\label{ssec-invrel-alg-ext-struct-alg}

Strictly following the ideas of Shepler~\cite{shepler}, we construct an algebra structure on $\Om^\chi$.
The first step is to define a product. For this, we use the polynomial $Q_\chi$ of definition~\ref{dfn-polynome-representation}
to bring back  the usual product of two elements of $\Om^\chi$ into $\Om^\chi$ (by Stanley theorem~\cite{stanley}, $Q_\chi$ is so that $S(V^*)=Q_\chi S(V^*)^G$).
Thus we are looking for divisibility conditions by $Q_\chi$ or more precisely by the non invertible part of $Q_\chi$ :
this is done in lemmas~\ref{lem-divisibilite-tsv} and~\ref{lem-divisibilite-om}. The wanted divisibility is obtained under
hypotheses on $\Bscr$ (hypotheses~\ref{hyp-bon-hyperplan} and~\ref{hyp-condition-etoile}).

\begin{hypothesis}\label{hyp-bon-hyperplan}
Let us assume that $\Bscr$ contains every $M$-bad hyperplane. Equivalently, every hyperplane contained in $\Gscr$ is $M$-good.
\end{hypothesis}

\subsubsection{Divisibility}

To begin with, let us extend the following result of divisibility~\cite[lemma 1]{gutkin} to the ring $T^{-1}S(V^*)$.

\begin{lemma}[Divisibility in ${T^{-1}S(V^*)}$]\label{lem-divisibilite-tsv}
Let $x \in T^{-1}S(V^*)$, $H \in \Hyp$ and $i \in \intnn{1}{e_H}$. Let us assume that $s_H x ={\zeta_H}^i x$.
Then $x$ is divisible by ${\al_H}^{e_H-i}$.
\end{lemma}

The lemma is interesting only for $H \in \Gscr$ since for $H \in \Bscr$, the linear form $\al_H$ is invertible in $T^{-1}S(V^*)$.

\begin{proofproof}
Since $T$ is $G$-stable, we can write $x = P/Q$ with $Q \in T^G$. Since $S(V^*)$ is an integral domain, we deduce that
$s_H P = {\zeta_H}^i P$. The lemma 1 of~\cite{gutkin} shows that $P$ is divisible by ${\al_H}^{e_H-i}$ and so is $x$. \smartqed
\end{proofproof}

We continue our study of divisibility by the $\al_H$. Let us consider the case of $\Om$.

\begin{lemma}[Divisibility in ${\Om}$]\label{lem-divisibilite-om}
Let us choose $\mu \in (\Om^p)^\chi$ and $H \in \Hyp$. We fix $(y_1, \ldots, y_r)$ a basis of $M^*$ so that $s_H(y_i)= {\zeta_H}^{j_i} y_i$
for all $i \in \intnn{1}{r}$ (see notation~\ref{notation-ji}). We write
$$\mu = \Sumu{1 \leq i_1 < \cdots < i_p \leq r} \mu_{i_1,\ldots, i_p} \ y_{i_1} \wedge \cdots \wedge y_{i_p}
    \qquad \ins{avec} \quad \mu_{i_1, \ldots, i_p} \in T^{-1}S(V^*).$$

For $H \in \Gscr$, we have
\begin{rn}{or\ }
\item[or\ ] $0 \leq j_{i_1} + \cdots + j_{i_p} \leq e_H - 1 - n_H(\chi) \qquad \ins{and}\qquad
            {\al_H}^{j_{i_1} + \cdots + j_{i_p} + n_H(\chi)} \mid \mu_{i_1, \ldots, i_p}$
\item[or\ ] $e_H-n_H(\chi) \leq j_{i_1} + \cdots + j_{i_p} \leq 2e_H - 2 - n_H(\chi) \quad\ins{and}\quad
            {\al_H}^{j_{i_1} + \cdots + j_{i_p} + n_H(\chi)-e_H} \mid \mu_{i_1, \ldots, i_p}$.
\end{rn}
\end{lemma}

\begin{proofproof}
Since the family $(y_{i_1} \wedge \cdots \wedge y_{i_p})_{1 \leq i_1 < \cdots < i_p \leq r}$ is a
$T^{-1}S(V^*)$-basis of $\Om^p$, we have

$$\mu_{i_1,\ldots, i_p}\ y_{i_1} \wedge \cdots \wedge y_{i_p} \in (\Om^p)^\chi.$$
Thus $s_H(\mu_{i_1,\ldots, i_p} \ y_{i_1} \wedge \cdots \wedge y_{i_p}) = {\zeta_H}^{-n_H(\chi)} \mu_{i_1,\ldots, i_p}\ y_{i_1} \wedge \cdots \wedge y_{i_p}.$
In addition
$$\begin{rcl} s_H(\mu_{i_1,\ldots, i_p} \ y_{i_1} \wedge \cdots \wedge y_{i_p}) &=&
    s_H(\mu_{i_1,\ldots, i_p}) \ s_H(y_{i_1}) \wedge \cdots \wedge s_H(y_{i_p})\\[0.5ex] &=&
    {\zeta_H}^{j_{i_1} + \cdots + j_{i_p}} s_H(\mu_{i_1,\ldots, i_p}) \ y_{i_1} \wedge\cdots\wedge y_{i_p}\,.\end{rcl}$$
We then deduce $s_H(\mu_{i_1,\ldots, i_p})  = {\zeta_H}^{-j_{i_1} - \cdots - j_{i_p} - n_H(\chi)} \mu_{i_1,\ldots, i_p}.$
The hypothesis~\ref{hyp-bon-hyperplan} on $\Bscr$ and $\Gscr$ tells us
$$\forall\,H \in \Gscr, \qquad 0 \leq n_H(M) = \Sum{j=1}{r} j_i \leq e_H-1\,.$$
Hence $2 - 2 e_H \leq {-j_{i_1} - \cdots - j_{i_p} - n_H(\chi)} \leq 0\,.$ There are two cases to distinguish :
\begin{rn}{or\ }
\item[or\ ] $1 - e_H \leq {-j_{i_1} - \cdots - j_{i_p} - n_H(\chi)} \leq 0$ and if we set
            $f_H = e_H -j_{i_1} - \cdots - j_{i_p} - n_H(\chi)$, we have
            $$s_H(\mu_{i_1,\ldots, i_p}) = {\zeta_H}^{f_H} \mu_{i_1,\ldots, i_p} \qquad \ins{with}\qquad 1 \leq f_H \leq e_H;$$
            \noindent the lemma~\ref{lem-divisibilite-tsv} ensures us that $\mu_{i_1,\ldots, i_p}$ is divisible by
            ${\al_H}^{j_{i_1} + \cdots + j_{i_p} + n_H(\chi)}$;

\item[or\ ] $2 -2 e_H \leq {-j_{i_1} - \cdots - j_{i_p} - n_H(\chi)} \leq -e_H$ and if
            $f_H = 2e_H -j_{i_1} - \cdots - j_{i_p} - n_H(\chi)$, we have
            $$s_H(\mu_{i_1,\ldots, i_p}) = {\zeta_H}^{f_H} \mu_{i_1,\ldots, i_p} \qquad \ins{with}\qquad 2 \leq f_H \leq e_H;$$
            \noindent the lemma~\ref{lem-divisibilite-tsv} ensures us that $\mu_{i_1,\ldots, i_p}$ is divisible by
            ${\al_H}^{j_{i_1} + \cdots + j_{i_p} + n_H(\chi)-e_H}$. \smartqed
\end{rn}
\end{proofproof}

To assure the divisibility of $Q_\chi$, we strengthen the hypothesis~\ref{hyp-bon-hyperplan} by the following one.

\begin{hypothesis}\label{hyp-condition-etoile} The subset $\Bscr$ verifies the hypothesis~\ref{hyp-bon-hyperplan}
and $\Bscr$ contains every hyperplane that is not $(M,\chi)$-acceptable. Equivalently, every hyperplane in $\Gscr$ is $(M,\chi)$-acceptable and $M$-good.
\end{hypothesis}

We then obtain the following result of divisibility by $Q_\chi$ which is a refinement of lemma 2 of~\cite{shepler}.

\begin{corollary}[Divisibility in ${\Om^\chi}$]\label{cor-divisibilite-omchi} Let us assume hypothesis~\ref{hyp-condition-etoile}.
For $\mu, \om \in \Om^\chi$, we have $Q_\chi \mid \mu \wedge \om$.
\end{corollary}

\begin{proofproof} We fix $H \in \Gscr$ and we consider the same basis $(y_1, \ldots, y_r)$ of $M^*$ of
lemma~\ref{lem-divisibilite-om}. When $I = \{i_1, \ldots, i_p\}$ is a subset of $\intnn{1}{r}$ with
$1 \leq i_1 < \cdots < i_p \leq r$, we set $y_I = y_{i_1} \wedge \cdots \wedge y_{i_p}$. Now, we can write
$$\mu = \Sumu{I \subset \intnn{1}{r}}\, \mu_I\; y_I \qquad \ins{and} \qquad \om =
      \Sumu{I \subset \intnn{1}{r}}\, \om_I\; y_I \qquad \ins{with} \qquad \mu_I, \om_I \in T^{-1}S(V^*).$$
Hence
$$\mu \wedge \om = \Sumu{I \cap J = \emptyset}\, \ep_{I,J}\mu_I\om_J\; y_{I\cup J} \qquad \ins{with}\qquad
      \ep_{I,J} \in \{\pm 1\}.$$

Now, let us use lemma~\ref{lem-divisibilite-om}. For this, we choose two subsets $I,J$ of $\intnn{1}{r}$ with $I \cap J = \emptyset$.
\begin{rs}
\item
If $0 \leq \Sumu{i \in I}{j_i} < e_H - n_H(\chi)$ or $0 \leq \Sumu{i \in J}{j_i} < e_H - n_H(\chi)$ then $\mu_I$
or $\om_J$ is divisible ${\al_H}^{n_H(\chi)}$ and then $\mu_I\om_J$ is divisible by ${\al_H}^{n_H(\chi)}$.

\item If not $\Sumu{i \in I}{j_i} \geq e_H -n_H(\chi)$ and $\Sumu{i \in J}{j_i} \geq e_H -n_H(\chi)$ with
$I \cap J = \emptyset$ which contradicts hypothesis~\ref{hyp-condition-etoile}.
\end{rs}

\noindent Hence the product $\mu \wedge \om$ is divisible by ${\al_H}^{n_H(\chi)}$ for all $H \in \Gscr$.
Since the family $(\al_H)_{H \in \Gscr}$ is constituted with elements prime to each other,
$\mu \wedge \om$ is divisible by
$$\Prodd{H \in \Gscr}{} {\al_H}^{n_H(\chi)}\,.$$
\noindent In addition, for $H \in \Bscr$, the element ${\al_H}^{n_H(\chi)}$ is invertible in $T^{-1}S(V^*)$,
we finally obtain the divisibility of $\mu \wedge \om$ by
$$\Prodd{H \in \Bscr} {\al_H}^{n_H(\chi)} \Prodd{H \in \Gscr}{} {\al_H}^{n_H(\chi)} = Q_\chi\,.$$
\end{proofproof}

\subsubsection{Algebra structure}

When hypothesis~\ref{hyp-condition-etoile} is assumed, the corollary~\ref{cor-divisibilite-omchi}
applies. For $\mu, \om \in \Om^{\chi}$, we can define the twisted product $\curlywedge$ by
$$\mu \curlywedge \om = {Q_\chi}^{-1}\mu \wedge \om \in \Om.$$
\noindent Actually, we have $\mu \curlywedge \om \in \Om^{\chi}$ and thus we define a law $\curlywedge$  on $\Om^\chi$
which gives to $\Om^\chi$ a structure of an associative $(T^G)^{-1}S(V^*)^G$-algebra with unit element $Q_\chi$.
Now, we have to show that $(\Om^{\chi}, \curlywedge)$ is an exterior $(T^G)^{-1}S(V^*)^G$-algebra.
For this, we study the structure constants of $(\Om^{\chi}, \curlywedge)$ and we will show that they are those of an exterior algebra.
To have more simple notation, we set $R= S(V^*)^G$.

\subsection{Exterior Algebra}\label{ssec-invrel-algebre-exterieure}

In the previous subsection, under the hypothesis~\ref{hyp-condition-etoile}, we have constructed  an algebra structure on $\Om^\chi$.
In this subsection, we are looking for its isomorphism class. The proof is divided in two stages :
first, we give a necessary and suffisant condition for the structure constants of $\Om^\chi$ to be those of an exterior algebra
(proposition~\ref{prop-cns-algebre-exterieure}); subsequently, we show that this condition is verified (proposition~\ref{prop-verif-cns-algebre-exterieure}).
To this perspective, we begin to generalize Stanley's theorem to the ring $T^{-1}S(V^*)$.

\begin{corollary}[Stanley's Theorem in ${T^{-1}S(V^*)}$]\label{cor-stanley-etendu} We have
\begin{equation}\label{eq-stanley}
(T^{-1}S(V^*))^{\chi} = (T^G)^{-1}S(V^*)^G Q_\chi \qetq (\Om^r)^{\chi} = (T^G)^{-1}S(V^*)^G Q_{\chi \cdot\det_M} \ins{vol}_M
\end{equation}
\noindent where $\ins{vol}_M$ is a non-zero element of $\La^r(M^*)$ once for all fixed.
\end{corollary}

\begin{proofproof} This is an easy consequence of the usual Stanley's theorem and of the fact that $k \ins{vol}_M = \La^r(M^*)$ is
a linear representation of $G$ with character $\det_{M^*}$. \smartqed
\end{proofproof}

\begin{notation}
Let $U$ be a $(T^{G})^{-1}S(V^*)$-module. For $u,v \in U$, we denote by
$u \doteq v$ if there exists $x \in ((T^{G})^{-1}R)^\times$ so that $xu =v$. In particular, $u$ an $v$
generate the same $(T^{G})^{-1}R$-submodule.
\end{notation}

\begin{proposition}[Necessary and sufficient condition]\label{prop-cns-algebre-exterieure}
Let us assume hypothesis~\ref{hyp-condition-etoile}. For every $\om_1, \ldots, \om_r \in (\Om^1)^\chi$,
the following propositions are equivalent :

\begin{ri}{ii}
\item for all $p \in \intnn{1}{r}$, the family
      $(\om_{i_1} \curlywedge \cdots \curlywedge \om_{i_p})_{1 \leq i_1 < \cdots < i_p \leq r}$ is
      a $(T^G)^{-1}R$-basis of $(\Om^p)^\chi$;

\item \begin{equation}\label{eq-produit}
\om_{1}\curlywedge\cdots\curlywedge\om_{r} \doteq Q_{\chi \cdot \det_M}\ins{vol}_M.
\end{equation}
\end{ri}
\end{proposition}

\begin{proofproof}
\donne{(i)}{(ii)}  This is an easy consequence of~(\ref{eq-stanley}).

\donne{(ii)}{(i)} We set $K = \ins{Frac}(T^{-1}S(V^*))$. Let us show that
$\Fscr = (\om_{i_1} \curlywedge \cdots \curlywedge \om_{i_p})_{i_1 < \cdots < i_p}$
is free over $K$. Since $\om_{i_1} \curlywedge\cdots\curlywedge \om_{i_p} = {Q_\chi}^{1-p}\ \om_{i_1} \wedge\cdots\wedge \om_{i_p}$,
it suffices to show that the family
$(\om_{i_1} \wedge \cdots \wedge \om_{i_p})_{1 \leq i_1 < \cdots < i_p \leq r}$ is free over $K$.
For that, let us consider the relation
$$\Sumu{1 \leq i_1 < \cdots < i_p \leq r} r_{i_1, \ldots, i_p}\, \om_{i_1} \wedge \cdots \wedge \om_{i_p} = 0
      \qquad \ins{with} \qquad r_{i_1, \ldots, i_p} \in K.$$

\noindent We fix $I = \{i_1, \ldots, i_p\} \subset \intnn{1}{r}$ with $1 \leq i_1 < \cdots < i_p \leq r$
and we set $I^{c} = \{i_{p+1}, \ldots , i_r\}$ the complementary of $I$ in $\intnn{1}{r}$. Multiplying the preceding relation
by $\om_{i_{p+1}} \wedge \cdots \wedge \om_{i_r}$, we obtain
$$r_{i_1, \ldots, i_p}\om_{i_1}\wedge\cdots\wedge \om_{i_p}\wedge\om_{i_{p+1}}\wedge\cdots\wedge\om_{i_r}=0.$$
Hence $0=r_{i_1, \ldots, i_p} {Q_\chi}^{r-1} \om_{1} \curlywedge \cdots \curlywedge \om_{r} \doteq
    r_{i_1, \ldots, i_p} {Q_\chi}^{r-1} Q_{\chi \cdot \det_M}\, \ins{vol}_M.$
We then deduce that $r_{i_1, \ldots, i_p} = 0$ and thus the $K$-freeness of the family
$(\om_{i_1} \wedge \cdots \wedge \om_{i_p})_{1 \leq i_1 < \cdots < i_p \leq r}$\,.

Finally, the family $\Fscr$ is a basis of $K$-vector space $K \otimes \La^p(M^*)$. Thus, if $\mu \in (\Om^p)^\chi$,
there exists $r_{i_1, \ldots, i_p} \in K$ so that
$$\mu =
    \Sumu{1 \leq i_1< \cdots < i_p \leq r} r_{i_1, \ldots, i_p}\, \om_{i_1} \curlywedge \cdots \curlywedge \om_{i_p}.$$

\noindent Let us fix again $I = \{i_1, \ldots, i_p\} \subset \intnn{1}{r}$ with $i_1 < \cdots < i_p$ and
set $I^{c} = \{i_{p+1}, \ldots , i_r\}$ its complementary. By multiplying the defining relation of $\mu$ by
$\om_{i_{p+1}} \curlywedge \cdots \curlywedge \om_{i_r}$, we obtain
$\mu \curlywedge \om_{i_{p+1}} \curlywedge \cdots \curlywedge \om_{i_r} \in (\Om^r)^\chi$. Thus, with~(\ref{eq-stanley}),
there exists $f \in (T^{G})^{-1}R$ so that
$$\mu \curlywedge \om_{i_{p+1}}\curlywedge\cdots\curlywedge\om_{i_r} = f Q_{\chi \cdot \det_M}\,\ins{vol}_M.$$
In addition,
$$\mu \curlywedge \om_{i_{p+1}} \curlywedge \cdots \curlywedge \om_{i_r} =
  \ep \,r_{i_1, \ldots, i_p} \om_{1} \curlywedge \cdots \curlywedge \om_{r} \doteq \ep \, r_{i_1, \ldots, i_p}
  Q_{\chi \cdot \det_M}\, \ins{vol}_M \quad \ins{avec} \quad \ep \in \{\pm 1\}.$$
Hence $r_{i_1, \ldots, i_p} \doteq \ep f \in (T^{G})^{-1}R.$
Therefore, the family $\Fscr$ is a $(T^{G})^{-1}R$-basis of $(\Om^p)^\chi$. \smartqed
\end{proofproof}

Now, the aim is to show that every $(T^{G})^{-1}R$-basis of $(\Om^1)^\chi$ verifies condition~\ref{eq-produit}.
For this, we construct a family $(\nu_i)_{1 \leq i \leq r}$ of $G$-invariants and a family
$(\mu_i)_{1 \leq i \leq r}$ of $\det_{M^*}$-invariants verifying respectively
$$\nu_1 \wedge \cdots \wedge \nu_r \doteq Q_{M}\, \ins{vol}_M \qquad \ins{and} \qquad
    \mu_1 \wedge \cdots \wedge \mu_r \doteq (Q_{M^*})^{r-1}\, \ins{vol}_M.$$

\begin{proposition}[Invariants in ${\Om^1}$]\label{prop-element-invariant}
There exists a family $(\nu_1, \ldots, \nu_r) \in ((\Om^1)^{G})^r$ and a family $(\mu_1, \ldots, \mu_r) \in ((\Om^1)^{\det_{M^*}})^r$ verifying
$$\nu_1 \wedge \cdots \wedge \nu_r \in k^\times Q_{M}\, \ins{vol}_M \qquad \ins{and} \qquad
        \mu_1 \wedge \cdots \wedge \mu_r \in k^\times (Q_{M^*})^{r-1}\, \ins{vol}_M.$$
\end{proposition}

\begin{proofproof}
The proof is a rephrasing of Gutkin's theorem~\cite{gutkin} using the notion of minimal matrix evolved
by Opdam in~\cite[definition 2.2 and proposition 2.4 $(ii)$]{opdam}.

For $C = (c_{ij})_{i,j} \in \mat{r}{S(V^*)}$, we denote by $g \cdot C$ the matrix $(g \,c_{ij})_{i,j}$.
Let us consider $C$ a $M$-minimal matrix. By definition, $C \in \mat{r}{S(V^*)}$ verifies
\begin{ri}{(iii)}
\item $g \cdot C = Cg_{M}$;
\item $\det C  \neq 0$;
\item $\deg \det C$ is minimal among the matrix verifying~$(i)$ and~$(ii)$.
\end{ri}

Let us choose $(y_i)_{1 \leq i \leq r}$ a basis of $M^*$ and define for $j \in \intnn{1}{r}$,
$$\nu_j = \Sum{i=1}{r} c_{ji} \otimes y_i.$$
We have $\nu_1 \wedge \cdots \wedge \nu_r \in k^\times \det(C) \ins{vol}_M.$
Although $\det C \in k^\times Q_M$ (let us see~\cite[proposition~2.4 $(iii)$]{opdam}).
It remains to show that $\nu_j$ is $G$-invariant for all $j \in \intnn{1}{r}$. But, $(i)$ gives
$$g\,c_{ji} = \Sum{k=1}{r} c_{jk}\, {g_M}_{ki} \qquad \ins{and} \qquad g\,y_i = \Sum{n=1}{r} {g_{M^*}}_{n i}\, y_n.$$
Hence
$$g \nu_j \ = \
    \Sum{i=1}{r}\pa{\Sum{k=1}{r} c_{jk}\, {g_M}_{ki} \otimes \Sum{n=1}{r} {g_{M^*}}_{n i}\,  y_n} \ = \
    \Sum{k=1}{r}\Sum{n=1}{r} \pa{\Sum{i=1}{r} {g_M}_{ki}{g_{M^*}}_{n i}} c_{jk} \otimes  y_n\,.$$

\noindent Since $\trans{g_{M^*}} = {g_M}^{-1}$, we have $\,\Sum{i=1}{r}\, {g_M}_{ki}{g_{M^*}}_{n i} = \de_{k n}$ and then
$$g \nu_j  = \Sum{k=1}{r} c_{jk} \otimes y_k = \nu_j.$$
\noindent Finally $\nu_j$ is $G$-invariant.

Now, let us consider $D$ a $M^*$-minimal matrix. By definition, $D \in \mat{r}{S(V^*)}$ verifies
\begin{ri}{(iii)}
\item $g \cdot D = Dg_{M^*}$;
\item $\det D  \neq 0$;
\item $\deg \det D$ is minimal among the matrix verifying~$(i)$ and~$(ii)$.
\end{ri}

We consider $\ins{Com\,} D = (e_{ij})_{i,j}$ the comatrix of $D$. Since the action of  $G$ on $S(V^*)$ is
compatible with the algebra structure, we have
$$g \cdot \ins{Com\,} D = \ins{Com\,} (g \cdot D) = \ins{Com\,} (Dg_{M^*}) = \ins{Com\,}D \ins{Com\,} g_{M^*}
    = \det(g_{M^*}) \ins{Com}(D)\,g_M.$$
We then define, for  $j \in \intnn{1}{r}$,
$$\mu_j \ =\ \Sum{i=1}{r} e_{ji} \otimes y_i\,.$$

We have $\mu_1 \wedge \cdots \wedge \mu_r \in k^\times \det(\ins{Com\,}D)\ins{vol}_M.$
Although $\det\ins{Com\,}D = (\det D)^{r-1}$ and $\det D \in k^\times Q_{M^*}$ (see~\cite[proposition 2.4 $(iii)$]{opdam}).
It remains to show that $\mu_j$ is $\det_{M^*}$-invariant for all $j \in \intnn{1}{r}$. But,
$$g\,e_{ji} = \det g_{M^*}\, \Sum{k=1}{r} e_{jk}\, {g_M}_{ki} \qquad \ins{and}\qquad
        g\,y_i = \Sum{n=1}{r} {g_{M^*}}_{n i}\, y_n.$$
\noindent Hence
$$g \mu_j\ =\
    \det g_{M^*} \Sum{i=1}{r}\pa{\Sum{k=1}{r} e_{jk}\, {g_M}_{ki} \otimes \Sum{n=1}{r} {g_{M^*}}_{n i}\, y_n}\ =\
    \det g_{M^*} \Sum{k=1}{r}\Sum{n=1}{r} \pa{\Sum{i=1}{r} {g_M}_{ki}{g_{M^*}}_{n i}  } e_{jk} \otimes y_n \,.$$
\noindent Since $\trans{g_{M^*}} = {g_M}^{-1}$, we have $\,\Sum{i=1}{r} {g_M}_{ki}{g_{M^*}}_{n i} = \de_{k n}$ and then
$$g \mu_j  = \det g_{M^*}\, \Sum{k=1}{r} e_{jk} \otimes y_k = \det g_{M^*} \mu_j.$$
\noindent Finally $\mu_j$ is $\det_{g_{M^*}}$-invariant. \smartqed
\end{proofproof}

The following lemma proposes polynomial relations : the aim (with an eye to proposition~\ref{prop-verif-cns-algebre-exterieure})
is to obtain formulas to determine when $Q_\chi Q_{M}(Q_{\chi \cdot \det_M})^{-1}$ and $Q_{\chi \cdot \det_M}Q_{M^*}(Q_\chi)^{-1}$ are prime to each other.

\begin{lemma}[Polynomial identities]\label{lem-egalite-polynomiale}
We define
$$\Gscr_0=\{H \in \Gscr,\quad n_H(M)=0\},\qquad\Gscr_+=\{H \in \Gscr,\quad n_H(M)\geq e_H - n_H(\chi)\}\,,$$
$$\Gscr_{\neq 0} = \Gscr \setminus \Gscr_0 \qqetqq \Gscr_- = \Gscr \setminus \Gscr_+.$$

We have
\begin{ri}{(iii)}
\item $$Q_{M^*}\ =\ \Prodd{H \in \Bscr} {\al_H}^{n_H(M^*)}
                          \Prodd{H \in \Gscr_{\neq 0}} {\al_H}^{e_H(r-n_{0,H}(M)) - n_H(M)}\,.$$
\item $$Q_{\chi \cdot \det_M}\ =\ \Prodd{H \in \Bscr} {\al_H}^{n_H(\chi \cdot \det_M)}
            \Prodd{H \in \Gscr_-} {\al_H}^{n_H(\chi) + n_H(M)}
            \Prodd{H \in \Gscr_+} {\al_H}^{n_H(\chi) + n_H(M)-e_H}\,.$$
\item $$Q_\chi Q_{M}(Q_{\chi \cdot \det_M})^{-1}\ =\
      \Prodd{H \in \Bscr} {\al_H}^{n_H(\chi)+ n_H(M)- {n_H(\chi \cdot \det_M)}}
      \Prodd{H \in \Gscr_+} {\al_H}^{e_H}\,.$$
\item $Q_{\chi \cdot \det_M}Q_{M^*}(Q_\chi)^{-1}\ =\ $
      $$\Prodd{H \in \Bscr}\!\! {\al_H}^{n_H(M^*) + n_H(\chi \cdot \det_M) - n_H(\chi)}\!\!\!\!
      \Prodd{H \in \Gscr_-\setminus \Gscr_0}\!\!\!\!\!\!{\al_H}^{e_H(r-n_{0,H}(M))}\!\!
      \Prodd{H \in \Gscr_+}\!\!{\al_H}^{e_H(r-1-n_{0,H}(M))}\,.$$
\end{ri}
\end{lemma}

\begin{proofproof} We have seen in the remark~\ref{rem-bon-excellent-hyperplan} that $n_H(M^*) = e_H(r-n_{0,H}(M)) -n_H(M)$
for every $H \in \Hyp$. Moreover, $n_H(M) = 0$ if and only if $n_H(M^*) = 0$ if and only if
$n_{0,H}(M) = r$. We then obtain~$(i)$.

Let $H \in \Gscr$. We have $0 \leq n_H(M) \leq e_H -1$ and then $n_H(M) = n_H(\det_M)$.
We conclude that $(\chi \cdot \det_{M})(s_H) = {\det (s_H)}^{- n_H(\chi) - n_H(M)}$.
Since $0 \leq n_H(\chi) + n_H(M) \leq 2 e_H -2$, we obtain
$n_H(\chi \cdot \det{}_{M}) = n_H(\chi) + n_H(M)\ \ins{ if }\ n_H(\chi) + n_H(M) \leq e_H-1$
and
$$n_H(\chi \cdot \det{}_{M}) = n_H(\chi) + n_H(M) - e_H \qquad \ins{if} \qquad n_H(\chi) + n_H(M) \geq e_H.$$
Identities~$(iii)$ and~$(iv)$ are easy consequences of~$(i)$ and~$(ii)$. \smartqed
\end{proofproof}

So that we can conclude on the algebra structure of $\Om^\chi$, we need to reinforce hypothesis.

\begin{hypothesis}\label{hyp-excellent-hyperplan} The subset $\Bscr$ contains all the hyperplanes that are not $(M,\chi)$-good that is
to say every element of $\Gscr$ that are $(M,\chi)$-good or equivalently $\Gscr_+ = \emptyset$.
\end{hypothesis}

\begin{hypothesis}\label{hyp-action-reflexion}
The subset $\Bscr$ contains all hyperplanes that are not $M$-excellents, or equivalently that $s_H$ acts on $M$ as identity or as a reflection for all $H \in \Gscr$.
\end{hypothesis}

\begin{remark}[Links between hypotheses]\label{rem-lien-hypothese} The remark~\ref{rem-bon-excellent-hyperplan} ensures that both hypotheses~\ref{hyp-excellent-hyperplan}
and~\ref{hyp-action-reflexion} are stronger than hypothesis~\ref{hyp-condition-etoile}.
In addition, under hypothesis~\ref{hyp-excellent-hyperplan}, the lemma~\ref{lem-egalite-polynomiale} shows
that $Q_\chi Q_{M} (Q_{\chi \cdot \det_M})^{-1}$ are invertible in $(T^{G})^{-1}R$. \smartqed
\end{remark}

\begin{proposition}[Checking of the necessary and sufficient condition]\label{prop-verif-cns-algebre-exterieure}
Let us assume that one of the two hypotheses~\ref{hyp-excellent-hyperplan} or~\ref{hyp-action-reflexion} are verified.
If $\om_1, \ldots, \om_r$ generate $(\Om^1)^{\chi}$ then
$$\om_1 \curlywedge \cdots \curlywedge \om_r \doteq Q_{\chi \cdot \det_M} \ins{vol}_M.$$
\end{proposition}

\begin{proofproof} From remark~\ref{rem-lien-hypothese}, the hypothesis~\ref{hyp-condition-etoile} is verified. Thus we can define the algebra structure on
$(T^{-1}S(V^*) \otimes \La(M^*))^{\chi}$. Since $\Om^\chi$ is stable by $\curlywedge$, we have $\om_1 \curlywedge \cdots \curlywedge \om_r \in (\Om^r)^\chi$.
The identity~(\ref{eq-stanley}) tells us
$$\exists f \in (T^{G})^{-1}R, \qquad
      \om_1 \curlywedge \cdots \curlywedge \om_r = f Q_{\chi \cdot \det_M} \ins{vol}_M.$$
Now, we have to prove that $f$ is invertible in $(T^{-1}S(V^*))^G$. But, actually it suffice to show that $f$ is invertible in $T^{-1}S(V^*)$.
Let us consider $(y_i)_{1 \leq i \leq r}$ a basis of $M^*$. We denote by $C \in \mat{r}{T^{-1}S(V^*)}$
the matrix of the family $(\om_i)_{1 \leq i \leq r}$ in the $T^{-1}S(V^*)$-basis $(1 \otimes y_i)_{1 \leq i \leq r}$
of $\Om^1$. We deduce the existence of $\la \in k^\times$ so that
$$\om_1 \wedge \cdots \wedge \om_r = \la \det C\, \ins{vol}_M \qquad \textrm{then} \qquad
        \la\det C =  f Q_{\chi \cdot \det_M} (Q_\chi)^{r-1}.$$

In addition, since $\nu_i$ is $G$-invariant, $Q_\chi \nu_i$ is $\chi$-invariant. We then deduce
that $Q_\chi \nu_i$ is a linear combination (with coefficients in $(T^{G})^{-1}R$\,) of the family $(\om_i)_{1 \leq i \leq r}$.
Thus we obtain a matrix $D \in \mat{r}{(T^{G})^{-1}R}$ so that
$$Q_\chi \nu_1 \wedge \cdots \wedge Q_\chi \nu_r  = \det D\, \om_1 \wedge \cdots \wedge \om_r =
      \la \det D \det C\, \ins{vol}_M = \la \det D\, f Q_{\chi \cdot \det_M} (Q_\chi)^{r-1} \ins{vol}_M$$

\noindent But $\nu_1 \wedge \cdots \wedge \nu_r \in k^\times Q_{M} \ins{vol}_M$, so we obtain
$Q_\chi \nu_1 \wedge \cdots \wedge Q_\chi \nu_r \in k^\times (Q_\chi)^r Q_{M} \ins{vol}_M$.
Hence
$$Q_\chi Q_{M} \in k^\times f Q_{\chi \cdot \det_M} \det D.$$
Finally $f \neq 0$ and $f \mid Q_\chi Q_{M} (Q_{\chi \cdot \det_M})^{-1}.$

Under hypothesis~\ref{hyp-excellent-hyperplan}, the remark~\ref{rem-lien-hypothese} shows that $f$ is invertible
and~(\ref{eq-produit}) is verified.

Now, let us assume hypothesis~\ref{hyp-action-reflexion}.
For $H \in \Gscr_+$, we have $n_{0,H}(M) = r-1$. Lemma~\ref{lem-egalite-polynomiale} shows that
$\pa{Q_{\chi \cdot \det_M} Q_{M^*} (Q_\chi)^{-1}}^{r-1}$ and $Q_\chi Q_{M}(Q_{\chi \cdot \det_M})^{-1}$ are prime to each other.
So that, we can conclude by showing that
$f \mid \pa{Q_{\chi \cdot \det_M} Q_{M^*} (Q_\chi)^{-1}}^{r-1}$.

Let $(\mu_i)_{1 \leq i \leq r}$ be the family of lemma~\ref{prop-element-invariant}.
Since $\mu_i$ is $\det_{M^*}$-invariant, $Q_{\chi \cdot \det_M} \mu_i \in (\Om^1)^\chi$.
Then $Q_{\chi \cdot \det_M} \mu_i$ is a linear combination (with coefficients in $(T^{G})^{-1}R$) of
$(\om_i)_{1 \leq i \leq r}$. By this way, we have constructed a matrix $D' \in \mat{r}{(T^{G})^{-1}R}$ so that
$$Q_{\chi \cdot \det_M} \mu_1 \wedge \cdots \wedge Q_{\chi \cdot \det_M} \mu_r =
      \det D'\, \om_1 \wedge \cdots \wedge \om_r = \det D' f Q_{\chi \cdot \det_M} (Q_\chi)^{r-1} \ins{vol}_M.$$

\noindent But $\mu_1 \wedge \cdots \wedge \mu_r \in k^\times (Q_{M^*})^{r-1}\; \ins{vol}_M$, so we obtain
$$Q_{\chi \cdot \det_M} \mu_1 \wedge \cdots \wedge Q_{\chi \cdot \det_M} \mu_r \in k^\times
  (Q_{\chi \cdot \det_M})^{r} (Q_{M^*})^{r-1}\,\ins{vol}_M\,.$$
Hence $(Q_{\chi \cdot \det_M}Q_{M^*})^{r-1} \in k^\times \det D'\, f (Q_\chi)^{r-1}.$
Finally $f \mid \pa{Q_{\chi \cdot \det_M} Q_{M^*} (Q_\chi)^{-1}}^{r-1}.$
Thus, $f$ divides $Q_\chi Q_{\det_M}(Q_{\chi \cdot \det_M})^{-1}$ and
$\pa{Q_{\chi \cdot \det_M} Q_{M^*} (Q_\chi)^{-1}}^{r-1}$ and identity~(\ref{eq-produit})
is verified under hypothesis~\ref{hyp-action-reflexion}. \smartqed
\end{proofproof}

\begin{theorem}[Exterior Algebra]\label{prop-conclusion-alg-exterieure} Let us assume that one of the two hypotheses~\ref{hyp-excellent-hyperplan}
or~\ref{hyp-action-reflexion} are verified. The $(T^{G})^{-1}R$-algebra $(\Om^{\chi}, \curlywedge)$ is an exterior algebra.
\end{theorem}

\begin{proofproof} From propositions~\ref{prop-cns-algebre-exterieure} and~\ref{prop-verif-cns-algebre-exterieure} and the remark~\ref{rem-lien-hypothese},
it suffices to show that $(\Om^1)^{\chi}$ can be generated by $r$ elements. Actually, we will show that $(\Om^1)^{\chi}$ is
a free module of rank $r$ over $(T^{G})^{-1}R$. By theorem B of Chevalley~\cite{chevalley}, we have
$$(S(V^*) \otimes M^*)^\chi = (S(V^*) \otimes M^* \otimes {k_\chi}^*)^G \otimes k_\chi = (S(V^*) \otimes (M \otimes k_\chi)^*)^G \otimes k_\chi\,.$$
\noindent
Thus we obtain
$$(S(V^*) \otimes M^*)^\chi = R \otimes (S_G \otimes (M \otimes k_\chi)^*)^G \otimes k_\chi$$
\noindent and $(S(V^*) \otimes M^*)^\chi$ is a free module of rank $\dim_k(\Hom{G}{S_G}{M \otimes k_\chi}) = \dim_k(M \otimes k_\chi) =r$.
By extending the scalar to $(T^G)^{-1} R$, we obtain that $(\Om^1)^\chi$ is free of rank $r$.  \smartqed
\end{proofproof}

\begin{remark}[Shepler, Orlik and Solomon]\label{rem-shepler-os} If every hyperplane of $\Hyp$ is $(M,\chi)$-good, we can choose $\Bscr = \emptyset$.
Similarly, if $s_H$ acts on $M$ as a reflection or as identity for all $H \in \Hyp$, we can choose $\Bscr = \emptyset$ and thus $T^{-1}S(V^*) = S(V^*)$.
We obtain the results of~\cite{beck} back and thus those of~\cite{orlik-solomon} and~\cite{shepler}. \smartqed
\end{remark}

\begin{remark}[When ${\Bscr = \Hyp}$] When $\Bscr = \Hyp$, the hypotheses~\ref{hyp-excellent-hyperplan} and~\ref{hyp-action-reflexion} are verified.
Thus $\Om^\chi$ is an $(T^{G})^{-1}R$-exterior algebra.  \smartqed
\end{remark}

\section{Consequences of the exterior algebra structure}\label{sec-invrel-consq-alg-exterieure}

In this section, we take an interest in the consequences of the structure of $\Om^\chi$ when $\Bscr$ is empty.
The first of these consequences is an equality between rational fonctions (corollary~\ref{cor-fraction-rationnelle}) generalizing
the one of Orlik and Solomon~\cite[equality 3.7]{orlik-solomon}.
In the subsection~\ref{ssec-consq-regulier}, we give various polynomial identities generalizing those
of~\cite{orlik-solomon}, \cite{lehrer-michel}, \cite{lehrer} and~\cite{bonnafe-lehrer-michel}. These identities leads to characterizations of the regularity of integers.

\begin{hypothesis} In this section, we assume that $\Bscr=\emptyset$.
Equivalently, we suppose that every hyperplane in $\Hyp$ is $(M,\chi)$-good or that $s_H$ acts on $M$ trivially or as a reflection for all $H \in \Hyp$.

Thus $(S(V^*) \otimes \La(M^*))^\chi$ is an $S(V^*)^G$-exterior algebra.
\end{hypothesis}

\subsection{Introduction and notations}\label{ssec-consq-intro-notation}

In this subsection, we introduce the objects studied next, in particular we set $\ga$ an element of the normalizer of $G$ in $\GL{V}$.
In addition, since the product of $\Om^\chi$ is a deformation of the usual product, we define a new degree which considers the deformation by $Q_\chi$
so that we obtain a bigraduation compatible with the algebra structure.

\begin{notation}[Bigraduation] Let us consider $S_n \subset S(V^*)$ the vector space of homogeneous polynomial functions with degree $n$.
For $p \in \intnn{0}{r}$, we set $\Om^p = S(V^*) \otimes \La^p(M^*)$ and $\Om_n^p = S_n \otimes \La^p(M^*)$. Thus, we have
$$\Om^\chi = \Somd{p=0}{r}\, (\Om^p)^\chi \qquad \ins{and} \qquad \Om^\chi = \Somdd{n \in \NN}\Somd{p=0}{r}\, (\Om_n^p)^\chi.$$
For $\om \in (\Om_n^p)^\chi$, we set $\deg(\om) = (n,p)$ and $\deg'(\om) = (n-\deg Q_\chi, p)$. If $\mu \in (\Om_{n'}^{p'})^\chi$ then
$$\om \curlywedge \mu \in (\Om_{n+n'-\deg Q_\chi}^{p+p'})^\chi$$
and $\deg'(\om) + \deg'(\mu) = \deg'(\om \curlywedge \mu)$.
\end{notation}

\begin{definition}[Fake degree, exponents and degrees] Let us remind what are the exponents of a representation of a reflection group.
We denote by $S_G$ the coinvariant ring of $G$ that is to say the quotient ring of $S(V^*)$ by the ideal generated by the polynomial invariant
functions vanishing at the origin. This is a graded $G$-module which is isomorphic to the regular representation (this is Chevalley's theorem~\cite{chevalley}).
We denote by $(S_G)_i$ the graded component of $S_G$ of degree $i$ and then define the fake degree of $M$ to be the polynomial
$$F_M(T) = \Sumu{i \in \NN} \langle (S_G)_i, M \rangle_G T^i \in \ZZ[T].$$
Since $\langle (S_G)_i, M \rangle_G$ are non negative integers, we can write $F_M(T) = T^{m_1(M)}+\cdots +T^{m_r(M)}$.
The integers $m_1(M), \ldots, m_r(M)$ are called the $M$-exponents.

By a theorem of Shephard and Todd~\cite{shephard-todd}, the ring $S(V^*)^G$ is generated by a family $(f_1,\ldots, f_\ell)$ of
$\ell$ homogeneous algebraically free generators. We denote by $d_i= \deg f_i$. The multiset $(d_1,\ldots, d_r)$  is well determined
and called the set of invariant degrees of $G$.
\end{definition}

The following lemma will be useful to extend the character $\chi$ of the group $\langle G,\ga\rangle$.

\begin{lemma}[Extension]\label{lem-prolongement}
Let $M,N$ and $P$ be three abelian groups and $\ph \colon M \ra N$, $\theta \colon M \ra P$ be two morphisms of groups.
We assume that $P$ is a divisible group. If $\ker \ph \subset \ker \theta$, there exists a morphism of groups $\wt{\theta} : N \ra P$
so that the following diagram commutes
$$\xymatrix{M \ar[r]^{\ph} \ar[d]^{\theta} & N \ar@{.>}[ld]^{\wt{\theta}}\\ P }$$
\end{lemma}

\begin{proofproof}
Since $\ker \ph \subset \ker \theta$, we can define a group homomorphism
$\theta_1 : M/\ker \ph \simeq \im \ph \ra P$ so that $\theta_1 \circ \ph = \theta$. Since $P$ is divisible,
we can extend $\theta_1$ in $\wt{\theta} : N \ra P$. Finally we obtain
$\wt{\theta} \circ \ph=\theta_1 \circ \ph = \theta$. \smartqed
\end{proofproof}

Let us introduce some notations and consider the normalizer $\Nscr$ of $G$ in $\GL{V}$.
We choose a semisimple element $\ga \in \Nscr$ (see~\cite{bonnafe-lehrer-michel}). We assume that $M$ is a $\langle G,\ga \rangle$-module and that $\ga$ acts
semisimply on $M$.
Furthermore we assume that the derived group $D$ of $\langle G,\ga \rangle$ verifies $D \subset \ker \chi$. By applying lemma~\ref{lem-prolongement} with
$M = G/D(G)$, $N = \langle G,\ga\rangle/D$, $P = \UU$ the groups of complex numbers with module $1$ and $\theta = \chi$,
we extend $\chi$ in a linear character of $\langle G,\ga \rangle$ (also denoted by $\chi$).
$$\xymatrix{G \ar[r] \ar[d] & \langle G,\ga\rangle \ar[d]\\ G/D(G) \ar[r]^{\ph} \ar[d]^\chi &
          \langle G,\ga\rangle /D \ar@{.>}[ld]^{\chi}\\ \UU}$$

We denote by $k_\chi$ the representation (of $\langle G,\ga \rangle$) with character $\chi$ over $k$ and we define
$M_\chi = M \otimes k_\chi$. So $M_\chi$ is an $\langle G,\ga \rangle$-module and, thanks to theorem B of Chevalley~\cite{chevalley},
we obtain an isomorphism of graded $G$-modules and of $R$-modules
$$(\Om^1)^\chi = (S(V^*) \otimes M^*)^\chi = (S(V^*) \otimes M^* \otimes {k_\chi}^*)^G \otimes k_\chi = R \otimes (S_G \otimes {M_\chi}^*)^G \otimes k_\chi\,.$$
Thus by definition of the $M_\chi$-exponents, we can choose an $R$-basis $\Gscr = (\om_1, \ldots, \om_r)$
of $(\Om^1)^\chi$ bihomogeneous with degree $\deg'(\om_i) =(m_i(M_\chi) - \deg(Q_\chi), 1)$. Moreover, the hypothesis
$D \subset \ker \chi$ ensure that $\ga$ stabilizes the vector space $N^\chi$ of $\chi$-invariants of~$N$,
for all $\langle G,\ga \rangle$-module $N$. Thus we obtain the isomorphism of graded $\langle \ga \rangle$-modules and of $R$-modules
$$(\Om^1)^\chi = (S(V^*) \otimes M^*)^\chi = (R \otimes S_G \otimes M^*)^\chi = R \otimes (S_G \otimes M^*)^\chi\,.$$
Finally, we can assume that the $\om_i$ are eigenvectors for $\ga$.
We denote by $\ep_{i,\ga, \chi}(M)$ the eigenvalue of $\ga$ associated to $\om_i$.
Both isomorphisms given above show that the multiset $(\ep_{i,\ga, \chi}(M),m_i(M_\chi))_i$
does not depend of the choice of the basis of $(\Om^1)^\chi$.

\begin{remark}[${m_i, \ep_i}$]
When $\chi = 1$ is the trivial character, we set $\ep_{i,\ga}(M) := \ep_{i,\ga, 1}(M)$.
Similarly, when $\ga = \id$, we set $\ep_{i,\chi}(M) := \ep_{i,\id, \chi}(M)$.
The family of $\ep_{i, \ga, \chi}(M)$ depends on the choice of the extension of $\chi$ to $\langle G, \ga \rangle$.

The family $\ep_{i, \ga}(V)$ can also be considered as the family of eigenvalues of $\ga$
so that the associated eigenvectors $(P_1, \ldots, P_\ell)$ are a family of
homogeneous and algebraically free generators of $R$ (see~\cite{bonnafe-lehrer-michel}). \smartqed
\end{remark}

\subsection{Rational Functions}\label{ssec-consq-fraction-rationnelle}

Here we follow the ideas of theorem~2.1 and equality~2.3 of~\cite{lehrer} and of the proposition~2.3 of~\cite{lehrer-michel}.

\begin{corollary}[Rational functions]\label{cor-fraction-rationnelle}
If $s_H$ acts on $M$ as the identity or as a reflection for all $H\in \Hyp$ or if $n_H(M) < e_H - n_H(\chi)$ for all $H \in \Hyp$ then

$$\frac{1}{|G|}\Sumu{g \in G} \overline{\chi(g)} \frac{\det(1 + (g\ga)_M Y)}{\det(1 - g\ga X)} =
      X^{\deg(Q_{\chi})} \frac{\Prod{i=1}{r} (1 + \ep_{i, \ga, \chi}(M)YX^{m_i(M_\chi) - \deg(Q_{\chi})})}
      {\Prod{i=1}{\ell} (1 - \ep_{i,\ga}(V) X^{d_i})}.$$
\end{corollary}

\begin{proofproof} The hypothesis $D \subset \ker \chi$ ensure that $\ga$ stabilizes $N^\chi$ the vector space of $\chi$-invariants of~$N$,
for all $\langle G,\ga \rangle$-module $N$. In particular, $\ga$ defines a bigraded endomorphism of $\Om^\chi$.
In order to show the equality, we compute the graded trace $P_{\Om^{\chi},\ga}(X,Y)$ of the endomorphism $\ga$
of $\Om^\chi$ in two different ways. By definition,
$$P_{\Om^{\chi},\ga}(X,Y) =\Sumu{n \in \NN}\,\Sum{p=0}{r} \tr\pa{\ga_{{\Om_n^p}^{\chi}}} X^nY^p\,.$$
Since $(\Om_n^p)^{\chi}$ is the $\chi$-isotypic component of $\Om_n^p$,
$$\frac{1}{|G|}\,\Sumu{g \in G}\,\overline{\chi(g)}\,g_{\Om_n^p}$$
\noindent is a projector on $(\Om_n^p)^{\chi}$. Hence
$$\tr\pa{\ga_{{\Om_n^p}^{\chi}}} =
      \frac{1}{|G|}\,\Sumu{g \in G}\,\overline{\chi(g)}\,\tr\pa{(g\ga)_{\Om_n^p}}.$$
Thus
$$P_{\Om^{\chi},\ga}(X,Y) = \frac{1}{|G|}\,\Sumu{g \in G}\, \overline{\chi(g)}
    \Sumu{n \in \NN}\,\Sum{p=0}{r} \tr\pa{(g\ga)_{\Om_n^p}} X^nY^p.$$
Finally, Molien's formulas give us
\begin{equation}\label{eq-alpha}
P_{\Om^{\chi},\ga}(X,Y) = \frac{1}{|G|}\,\Sumu{g \in G} \overline{\chi(g)}\,
      \frac{\det(1 + (g\ga)_M Y)}{\det(1 - g\ga X)}\,.
\end{equation}

In addition, propositions~\ref{prop-cns-algebre-exterieure} and~\ref{prop-verif-cns-algebre-exterieure} show
that $\Om^{\chi} = R \otimes \curlywedge ((\Om^1)^\chi)$ where $\curlywedge ((\Om^1)^\chi)$ is the $k$-algebra
(for $\curlywedge$) generated by $(\Om^1)^\chi$. Since the product $\curlywedge$ is compatible with $\deg'$
and since the degree of the unit element $e = Q_\chi$ for $\curlywedge$ is $\deg'(e) = (0,0)$, we obtain

$$P_{\curlywedge ((\Om^1)^\chi), \ga} (X,Y) =
          X^{\deg(Q_\chi)} \Prod{i=1}{r} (1 + \ep_{i,\ga,\chi}(M)Y X^{m_i(M_\chi) - \deg(Q_{\chi})}).$$
Moreover,
$$P_{R, \ga} (X) =\Prod{i=1}{\ell} (1 - \ep_{i, \ga}(V)X^{d_i})^{-1},$$
Hence
\begin{equation}\label{eq-beta}
P_{\Om^{\chi}, \ga}(X,Y) = X^{\deg(Q_\chi)}\,
            \frac{\Prod{i=1}{r} (1 + \ep_{i,\ga, \chi}(M)Y X^{m_i(M_\chi) - \deg(Q_{\chi})})}
            {\Prod{i=1}{\ell} (1 - \ep_{i,\ga}(V)X^{d_i})}.
\end{equation}
\noindent The equalities $\ref{eq-alpha}$ and $\ref{eq-beta}$ give the result. \smartqed
\end{proofproof}

\begin{remark}[${m_i(M_\chi)}$]
If $n_H(M) < e_H - n_H(\chi)$ for all $H \in \Hyp$, then the multisets
$$\paa{m_1(M) + \deg(Q_\chi), \ldots, m_r(M) + \deg(Q_\chi)} \qquad \ins{and} \qquad \paa{m_1(M_\chi), \ldots, m_r(M_\chi)}$$
are the same. Indeed, by following the proof of proposition~\ref{prop-verif-cns-algebre-exterieure},
we notice that, under our hypothesis, the family $(Q_\chi \nu_i)_{1 \leq i \leq r}$ is a basis of $(\Om^1)^\chi$.
But the properties of minimal matrices allow us to choose $\nu_i$ bihomogeneous with degree $(m_i(M),1)$. \smartqed
\end{remark}

\subsection{Regular Integers}\label{ssec-consq-regulier}

Similarly to the article of Lehrer and Michel~\cite{lehrer-michel} and the article of Lehrer~\cite{lehrer},
let us see apply identity~\ref{cor-fraction-rationnelle} to the representations  $V^\si$ and ${V^{*}}^{\si}$ where $\si \in \gal{\overline{\QQ}}{\QQ}$.
Let $d\in \NN$ and $\xi$ be a primitive $d^{\ins{th}}$ root of unity; we then define
$$A_\ga(d)=\{i \in \intnn{1}{\ell}, \quad \ep_{i,\ga}(V)\xi^{-d_i}=1\} \qquad \ins{and}\qquad a_\ga(d) = |A_\ga(d)|,$$
and for $\si \in \galq$,
$$r_i(\si, \chi) =\deg(Q_{\chi}) - m_i({V^\si}_\chi),
\quad r_i^*(\si, \chi) =\deg(Q_{\chi}) - m_i({{V^*}^\si}_\chi),$$
$$B_{\si,\ga}(d, \chi) = \{j \in \intnn{1}{\ell}, \quad
    \ep_{j,\ga,\chi}(V^{\si})\xi^{-\si}\xi^{r_j(\si, \chi)}=1\} \quad \ins{and}\quad b_{\si,\ga}(d, \chi) = |B_{\si,\ga}(d, \chi)|;$$
and
$$B^*_{\si,\ga}(d, \chi) =
            \{j \in \intnn{1}{\ell}, \quad \ep_{j,\ga,\chi}({V^*}^{\si})\xi^{\si}\xi^{r_j^*(\si, \chi)}=1\}
            \quad \ins{and} \quad b^*_{\si,\ga}(d, \chi) = |B^*_{\si,\ga}(d, \chi)|.$$
\noindent For $h \in \End{\CC}{V}$, we denote by $\det'(h)$ the product of non-zero eigenvalues of $h$,
we denote also by $V(h,\xi)= \ker(h-\xi \id)$ the eigenspace of $h$ associated to the eigenvalue $\xi$ and we set $d(h,\xi)= \dim(V(h,\xi))$.

\begin{theorem}\label{th-formule-pw}
We have $a_\ga(d) \leq \bsdchigamma$ and the following identity in $\CC[T]$

\noindent $\xi^{\deg(Q_\chi)}\Sumu{g \in G}\overline{\chi(g)}T^{d(g\ga,\xi)}(\det'(1-\xi^{-1}g\ga))^{\si-1}=$

\noindent $\begin{accolade} \Prodd{j \in \Bsdchigamma}\!\!\!\!\!\! (T-r_j(\si, \chi))
      \Prodd{j \notin \Bsdchigamma}\!\!\!\!\!\! (1-\ep_j\xi^{r_j(\si, \chi) -\si}) \Prodd{j \notin A_\ga(d)}
      \frac{d_j}{1-\ep'_j\xi^{-d_j}}\!\! & \textrm{ if } a_\ga(d)=\bsdchigamma,\\[3ex]  0 & \textrm{otherwise}\\
\end{accolade}$

\noindent where $\ep_i = \ep_{i,\ga, \chi}(V^\si)$ and $\ep'_i = \ep_{i,\ga}(V)$.

\vskip0.3cm
We have $a_\ga(d) \leq b^*_{\si,\ga}(d, \chi)$ and the following identity in $\CC[T]$

\noindent $(-1)^\ell\xi^{\deg(Q_\chi)+ \ell\si}\Sumu{g \in G}\overline{\chi(g)}(-T)^{d(g\ga,\xi)}
    (\det'(1-\xi^{-1}g\ga))^{\si-1}\det(g\ga)^{-\si}=$

\noindent $\begin{accolade} \Prodd{j \in \Bsdchigammastar}{}\!\!\!\!\! (T-r_j^*(\si,\chi))
    \Prodd{j \notin \Bsdchigammastar}\!\!\!\!\!\!\! (1-\ep_j\xi^{r_j^*(\si, \chi) +\si})\!\!
    \Prodd{j \notin A_\ga(d)} \frac{d_j}{1-\ep'_j\xi^{-d_j}}\!\!\!
    & \textrm{ if } a_\ga(d)=\bsdchigammastar,\\[3ex]  0 & \textrm{otherwise}\\
\end{accolade}$
\vskip0.2cm
\noindent where $\ep_i = \ep_{i,\ga, \chi}({V^*}^\si)$ and $\ep'_i = \ep_{i,\ga}(V)$.
\end{theorem}

\begin{proofproof}
For every reflection $s \in G$, $s_{V^\si}$ is still a reflection.
Thus we can apply corollary~\ref{cor-fraction-rationnelle} to the $G$-module $V^\si$. For $g \in G$,
we denote by $\la_1(g\ga), \ldots, \la_\ell(g\ga)$ the eigenvalues of $g\ga$ acting on $V$. We have
$$\frac{1}{|G|}\Sumu{g \in G}\overline{\chi(g)}\Prod{i=1}{\ell} \frac{(1+Y(\la_i(g\ga))^{\si})}{(1-X\la_i(g\ga))} =
X^{\deg(Q_{\chi})} \Prod{i=1}{\ell} \frac{(1 + \ep_iYX^{-r_i(\si, \chi)})}{(1 -\ep'_iX^{d_i})}.$$

\noindent We switch the indeterminate with $Y = \xi^{-\si}(T(1-\xi X)-1)$.

Let us begin with the left side. It becomes
$$\frac{1}{|G|}\Sumu{g \in G}\overline{\chi(g)}\Prod{i=1}{\ell}
      \frac{1-(\la_i(g\ga)\xi^{-1})^{\si}(1-T(1-\xi X))}{1-X\la_i(g\ga)}.$$
In each term of the sum, we discriminate the eigenvalues of $g\ga$ between those equal to $\xi$ and the others.
We obtain in $\CC(T,X)$
$$\frac{1}{|G|}\Sumu{g \in G}\overline{\chi(g)} \pa{\Prodd{\{i \,|\, \la_i = \xi\}} \!\!T
      \Prodd{\{i \mid \la_i \neq \xi \}} \!\! \frac{1-(\la_i(g\ga)\xi^{-1})^{\si}(1-T(1-\xi X))}{1-X\la_i(g\ga)}}.$$
So $\xi^{-1}$ is not a pole of this rational function with respect to $X$ and evaluating at $X=\xi^{-1}$, we obtain
$$\frac{1}{|G|}\Sumu{g \in G}\overline{\chi(g)}\,T^{d(g\ga,\xi)}\pa{\Prodd{\{i | \la_i \neq \xi\}}
      \frac{1-(\xi^{-1}\la_i(g\ga))^{\si}}{1-\xi^{-1}\la_i(g\ga)}} =\frac{1}{|G|} \Sumu{g \in G}\overline{\chi(g)}\,
      T^{d(g\ga,\xi)}\,(\det'(1-\xi^{-1}g\ga))^{\si-1}.$$

Now, let us consider the right side. After switching the indeterminate, it becomes
$$X^{\deg(Q_\chi)}\Prod{i=1}{\ell}\frac{1-\ep_i\xi^{-\si}(1-T(1-\xi X))X^{-r_i(\si, \chi)}}{1-\ep'_iX^{d_i}}.$$

\noindent Let us count the multiplicity of $\xi^{-1}$ as a root of the numerator and of the denominator of this rational function
with respect to $X$. For the denominator, $\xi^{-1}$ is a root of $1-\ep'_iX^{d_i}$ if and only if $i \in A_\ga(d)$. Moreover
this root is simple. So $\xi^{-1}$ is a root of order $a_\ga(d)$ of the denominator. For the numerator,
$$1-\ep_i\xi^{-\si}(1-T(1-\xi X))X^{-r_i(\si,\chi)}$$
is zero for $X=\xi^{-1}$ if and only if $i \in \Bsdchi$. Moreover, when differentiating with respect to $X$, we obtain
$$-\ep_i\xi^{-\si}\pa{-r_i(\si, \chi)(1-T(1-\xi X))X^{-r_i(\si, \chi)-1}+T\xi X^{-r_i(\si, \chi)}}$$
which is nonzero at $X=\xi^{-1}$. So $\xi^{-1}$ is a root of order $\bsdchigamma$ of the numerator.

Since $\xi^{-1}$ is not a pole of the left side, we obtain $a_\ga(d)\leq \bsdchigamma.$
Moreover, we deduce that the right side is zero if $a_\ga(d) < \bsdchigamma$.

Now, let us assume that $a_\ga(d)= \bsdchigamma$. If $i \in \Bsdchigamma$ then $\ep_i\xi^{-\si}=\xi^{-r_i(\si, \chi)}$
and
$$\begin{array}{r @{\ = \ } l} 1-\ep_i\xi^{-\si}(1-T(1-\xi X))X^{-r_i(\si, \chi)}
    & 1-(1-T(1-\xi X))(\xi X)^{-r_i(\si, \chi)} \\[1.5mm]
    & 1-(\xi X)^{-r_i(\si, \chi)} + T(1-\xi X)(\xi X)^{-r_i(\si,\chi)}\\[0.5mm]
    & (1-\xi X)\left(T(\xi X)^{-r_i(\si, \chi)} + \Sum{k = 0}{-r_i(\si, \chi)-1}(\xi X)^{k}\right).\end{array}$$

\noindent For $j \in A_\ga(d)$, we have $\ep'_j = \xi^{d_j}$ and so
$$\begin{array}{r @{\ = \ } l} 1-\ep'_j X^{d_j} & 1 - \xi^{d_j} X^{d_j}  \\
      & (1-\xi X)\ \Sum{k = 0}{d_j-1}(\xi X)^{k}. \end{array}$$
As a consequence, for $j \in A_\ga(d)$ and $i \in \Bsdchigamma$, we obtain
$$\frac{1-\ep_i\xi^{-\si}(1-T(1-\xi X))X^{-r_i(\si, \chi)}} {1-X^{d_j}}=\frac{T(\xi X)^{-r_i(\si, \chi)} +
    \Sum{k = 0}{-r_i(\si, \chi)-1}(\xi X)^{k} }{\Sum{k = 0}{d_j-1} (\xi X)^{k}}.$$

\noindent Evaluating at $X=\xi^{-1}$, we obtain $\frac{T -r_i(\si, \chi)}{d_j}$.
Finally, by choosing for each factor of the numerator whose index is in $\Bsdchigamma$,
one of the factor of the denominator whose index is in $A_\ga(d)$ (this is possible since $a_\ga(d)=\bsdchigamma$),
we obtain, after evaluating at $X=\xi^{-1}$,

$$\xi^{-\deg(Q_\chi)}\frac{\Prodd{j \in \Bsdchigamma}(T -r_j(\si, \chi)) \Prodd{j \notin \Bsdchigamma}
    (1 -\ep_j\xi^{r_j(\si, \chi)-\si})} {\Prodd{j \in A_\ga(d)} d_j\Prodd{j \notin A_\ga(d)} (1 -\ep'_j\xi^{-d_j})}.$$
\noindent The relation $|G| = \Prod{i=1}{\ell} d_i$ give us the identity.

\vskip0.3cm
For the second identity, we apply the corollary~\ref{cor-fraction-rationnelle} to $V{^{*}}^{\si} = {V^{\si}}^{*}$ on which $s_H$ acts
as a reflection for all $H \in \Hyp$. We switch the indeterminate in $Y = \xi^\si(T(1-\xi X)-1)$ and simplify with $(1-z^{-1})(1-z)^{-1} = - z^{-1}$. \smartqed
\end{proofproof}

\subsubsection{When ${\ga}$ is trivial}
We are interested in the case where $\ga = \id$. To simplify the notations, we set
$$B^*(d,\chi) := B^*_{\id,\id}(d,\chi) = \paa{j \in \intnn{1}{\ell},\quad d \mid 1 + r_j^*(\id,\chi)} \quad \ins{and} \quad
      b^*(d,\chi) = |B^*(d,\chi)|;$$
$$B(d,\chi) := B_{\id,\id}(d,\chi) = \paa{j \in \intnn{1}{\ell},\quad d \mid 1 - r_j(\id,\chi)} \qquad \ins{and} \qquad
      b(d,\chi) = |B(d,\chi)|;$$
and finally $A(d) := A_\id(d) = \paa{j \in \intnn{1}{\ell},\quad d \mid d_j}$ and $a(d) = |A(d)|$.

Let us remind that $d$ is said to be a regular integer if one of the $V(g, \xi)$ meets the complementary of the hyperplanes of~$\Hyp$.
The following corollary generalizes the results of~\cite{lehrer-michel} and the one of~\cite{lehrer}.

\begin{corollary}[Consequences and Exceptional Case]\label{cor-consequence-formule} We obtain the following formulas
\begin{ri}{viii}
\item $\Sumu{g \in G}\overline{\chi(g)} T^{d(g,1)}(\det'(1-g))^{\si-1} = \Prod{j=1}{\ell}(T - r_j(\si, \chi))$.

\item $\xi^{\deg(Q_\chi)}\Sumu{g \in G} \overline{\chi(g)}T^{d(g,\xi)}=$

      $\begin{accolade} \Prodd{j \in B(d,\chi)} (T-r_j(\id, \chi)) \Prodd{j \notin B(d,\chi)} (1-\xi^{r_j(\id,\chi)-1})
      \Prodd{j \notin A(d)} \frac{d_j}{1-\xi^{-d_j}},
      & \textrm{if } a(d)=b(d,\chi),\\[3ex]
      0 & \textrm{otherwise}.\end{accolade}$

\item $\Sumu{g \in G}\chi(g)T^{d(g,1)} = \Prodd{j=1}^{\ell}(T- r_j(\id,\chi))$.

\item We have  $a(d) \leq b^*(d,\chi)$ and

\noindent $(-1)^{\ell} \xi^{\ell + \deg(Q_\chi)}\Sumu{g \in G}(-T)^{d(g,\xi)}(\chi \cdot \det)(g^{-1})=$

\noindent $\begin{accolade} \Prodd{j \in B^*(d,\chi)}\!\! (T - r_j^*(\id,\chi)) \Prodd{j \notin B^*(d,\chi)}\!\!
    (1-\xi^{r_j^*(\id,\chi) + 1}) \Prodd{j \notin A(d)} \frac{d_j}{1-\xi^{-d_j}},
    &\ \ \textrm{if } a(d)=b^*(d,\chi),\\[3ex]
    0 &\ \ \textrm{otherwise}.\end{accolade}$

\item $\Sumu{g \in G}T^{d(g,1)}(\chi \cdot \det)(g)= \Prod{j=1}{\ell}(T + r_j^*(\id,\chi))$.

\item The multisets $\paa{\,-r_1^*(\id, \chi), \ldots, - r_\ell^*(\id, \chi)\,}$ and
      $\paa{\,r_1(\id, \chi \cdot \det), \ldots, r_\ell(\id, \chi \cdot \det)\,}$ are the same and
      $b^*(d,\chi) = b(d,\chi \cdot \det)$.

\item If $d$ is regular, then $a(d) = b_\si(d, \chi)$ for every $\si \in \gal{\overline{\QQ}}{\QQ}$ and every one dimensional character $\chi$.

\item If for all $H \in \Hyp$, the restriction of $\chi \cdot \det$ to $G_H$ is non trivial, then $d$ is a regular integer
      if and only if $a(d) = b(d,\chi)$.
\end{ri}
\end{corollary}

\begin{proofproof}
\begin{ri}{$viii$}
\item This is theorem~\ref{th-formule-pw} with $d=1$ and so $\xi=1$. We have
      $$A(1)=B_{\si,\id}(1,\chi)=\intnn{1}{\ell} \qqetqq a(1)=b_{\si,\id}(1,\chi).$$

\item This is theorem~\ref{th-formule-pw} with $\si=\id$. This Lehrer's identity~2.1~\cite{lehrer}.

\item Since $\overline{\chi(g)} = \chi(g^{-1})$ and $d(g,1) = d(g^{-1},1)$, this is $(i)$ for
      $\si=\ins{id}$ or $(ii)$ with $d=1$. This is Lehrer'identity 3.2~\cite{lehrer}.

\item This is the second identity of theorem~\ref{th-formule-pw} for $\si = \id$.

\item This is $(iv)$ with $d=1=\xi$ with the remark that $d(g,1)=d(g^{-1},1)$.

\item This one is obtained by comparing $(iii)$ applied to $\chi \cdot \det$ with $(v)$.

\item Theorem 3.4 of Springer~\cite{springer-article} shows us that the degree of the polynomial of the left side in
      theorem~\ref{th-formule-pw} is at most $a(d)$. Let us compute the coefficient of $T^{a(d)}$. Since $d$ is regular,
      the $g \in G$ verifying $d(g,\xi) = a(d)$ are a single conjugacy class. Thus, for every $g$ so that $d(g,\xi) = a(d)$,
      the value of $\overline{\chi(g)}\det'(1-\xi^{-1}g)^{\si - 1}$ does not depend on $g$.
      The coefficient of $T^{a(d)}$ is non zero and so $a(d) = b_\si(d,\chi)$.

\item This is $(vi)$ coupled with corollary 3.9 of Lehrer~\cite{lehrer} applied to the linear character $\chi\cdot\det$. \smartqed
\end{ri}
\end{proofproof}

\subsubsection{When ${\ga}$ is not necessarily trivial}
The case $\chi = 1$ is done in~\cite{bonnafe-lehrer-michel}. Let us remind that $d$ is $\ga$-regular
if one of the eigenspaces $V(g\ga, \xi)$ meets the complementary of the hyperplanes~$\Hyp$.
As a matter of simplification, we set $B_\ga^*(d,\chi) := B^*_{\id,\ga}(d,\chi)$, $b_\ga^*(d,\chi) := |B_\ga^*(d,\chi)|$ and
$$B_\ga(d,\chi) := B_{\id,\ga}(d,\chi)\qetq b_\ga(d,\chi) :=|B_\ga(d,\chi)|.$$

\begin{corollary}[Consequences and Exceptionnal Cases]\label{cor-consequence-formule-nouveau} We obtain the following formulas
\begin{ri}{iii}

\item $\xi^{\deg(Q_\chi)}\Sumu{g \in G} \overline{\chi(g)}T^{d(g\ga,\xi)}=$

      \noindent $\begin{accolade} \Prodd{j \in B_\ga(d,\chi)}\!\!\! (T-r_j(\id, \chi))\!\!
      \Prodd{j \notin B_\ga(d,\chi)}\!\!\! (1-\ep_j\xi^{r_j(\id,\chi)-1})
      \Prodd{j \notin A_\ga(d)}\! \frac{d_j}{1-\ep'_j\xi^{-d_j}},\!\!
          & \textrm{if } a_\ga(d)=b_\ga(d,\chi),\\[3ex]
        0 & \textrm{otherwise}.
      \end{accolade}$

      \noindent when $\ep_i = \ep_{i,\ga, \chi}(V)$ and $\ep'_i = \ep_{i,\ga}(V)$.

\item $(-1)^\ell\xi^{\deg(Q_\chi)+ \ell} \det(\ga^{-1})\Sumu{g \in G} (\chi\cdot\det)(g^{-1})(-T)^{d(g\ga,\xi)}=$

      \noindent $\begin{accolade} \Prodd{j \in B^*_\ga(d, \chi)}{}\!\!\!\! (T-r_j^*(\id,\chi))\!\!
      \Prodd{j \notin B^*_\ga(d, \chi)}\!\!\!\! (1-\ep_j\xi^{r_j^*(\id, \chi) +1})
      \Prodd{j \notin A_\ga(d)} \frac{d_j}{1-\ep'_j\xi^{-d_j}}\!\!
      & \textrm{ if } a_\ga(d)=b^*_\ga(d, \chi),\\[3ex]  0 & \textrm{otherwise}\\
      \end{accolade}$

      \noindent when $\ep_i = \ep_{i,\ga, \chi}({V^*})$ and $\ep'_i = \ep_{i,\ga}(V)$.

\item The two multisets $\paa{\,-r_i^*(\id, \chi), \quad i \in B^*_\ga(d, \chi)}$ and
      $\paa{\,r_i(\id, \chi \cdot \det), \quad i \in B_\ga(d, \chi)}$ are the same and $b_\ga^*(d,\chi) = b_\ga(d,\chi \cdot \det)$.

\item If $d$ is $\ga$-regular, then $a_\ga(d) = b_{\si,\ga}(d, \chi)= b^*_{\si,\gamma}(d,\chi)$ for every linear character $\chi$ and every
        $\si \in \gal{\overline{\QQ}}{\QQ}$.

\item If for all $H \in \Hyp$, the restriction of $\chi$ to $G_H$ is non-trivial, then $d$ is $\ga$-regular if and only if $a_\ga(d) = b_\ga(d, \chi)$.

\item If for all $H \in \Hyp$, the restriction of $\chi \cdot \det$ to $G_H$ is non-trivial, then $d$ is
      $\ga$-regular if and only if $a_\ga(d) = b^*_\ga(d, \chi)$.
\end{ri}
\end{corollary}

\begin{proofproof}
\begin{ri}{$iii$}
\item This is theorem~\ref{th-formule-pw} with $\si=\id$.

\item This is the second identity of theorem~\ref{th-formule-pw} with $\si=\id$.

\item Let us compare the roots of $(i)$ applied to $\chi \cdot \det$ and those of $(ii)$.

\item The theorem~3.4 of Springer~\cite{springer-article} show us that the degree of polynomial of the left side
      in theorem~\ref{th-formule-pw} is at most $a_\ga(d)$. Let us compute the coefficient of $T^{a_\ga(d)}$.
      Since $d$ is regular, the $g \in G$ so that $d(g\ga,\xi) = a_\ga(d)$ are a single conjugacy class. Thus the value
      $\overline{\chi(g)}\det'(1-\xi^{-1}g)^{\si - 1}$ does not depend on $g$ when $g$ verifies $d(g\ga,\xi) = a_\ga(d)$.
      The coefficient of $T^{a_\ga(d)}$ is non-zero and so $a_\ga(d) = b_{\si,\ga}(d,\chi)$.

\item By $(iv)$, it suffices to show that if $a_\ga(d) = b_\ga(d, \chi)$ then $d$ is $\ga$-regular.
      By $(i)$, the coefficient of $T^{a_\ga(d)} = T^{b_\ga(d, \chi)}$ is non-zero. Thus, thanks to Springer theorem,
      $\sum_{g \in C} \chi(g)$ is a factor of the coefficient of $T^{a_\ga(d)}$ where
      $$C = \{g \in G,\qquad \forall\,x \in V(h\ga,\xi), \quad gx=x\} \qquad \ins{with}\quad d(h\ga, \xi) = a_\ga(d).$$
      If $C$ is not the trivial group, then $C$ contains one of the $G_H$ (this is Steinberg's theorem) and since de restriction of $\chi$ to $G_H$ is non-trivial,
      we have $\sum_{g \in C} \chi(g) = 0$. Finally we obtain a contradiction and $C=1$ which means exactly that $d$ is $\ga$-regular.

\item This is $(iii)$ and $(v)$. \smartqed
\end{ri}
\end{proofproof}

\section{Types of hyperplanes}\label{sec-type-hyperplane}

In definition~\ref{dfn-type-hyperplan}, we define various types of hyperplanes. In this section, we study these types of hyperplane for some examples of reflection groups,
namely the symmetric group, the wreath product $G(d,1,n)$, the imprimitive groups of rank $2$ that is $G(de,e,2)$ and some exceptional case $G_4$, $G_5$ and $G_{24}$
(named after the classification of~\cite{shephard-todd}). The details of the computations can be found in~\cite{beck-these}.

\subsection{The symmetric group}

The symmetric groups $\Sfr_n$ acts faithfully as a reflection group over $\CC^{n}/\langle (1,\ldots, 1)\rangle$ by permuting the coordinates.
The reflections are the transpositions. They are of order 2 and conjugate to each other. Hence there is a unique conjugacy class of hyperplane.
The linear character of $\Sfr_n$ are the trivial one (denoted by $1$) and the sign character (denoted by $\ep$).
In addition, the representation of $\Sfr_n$ are described by the partition of $n$ (see for example~\cite{okounkov}).

\begin{proposition}[Symmetric group] Let $H$ be an hyperplane of the reflection group $\Sfr_n$ and $\rho$ be an irreducible representation of $\Sfr_n$.

The hyperplane $H$ is $\rho$-excellent, $\rho$-good, $(\rho,1)$-good, $(\rho,\ep)$-acceptable if and only if $\rho = 1$ or $\rho =\ep$ or $\rho$ is
the reflection representation or $n=4$ and $\rho$ is associated to the partition $(2,2)$.

The hyperplane $H$ is $(\rho,\ep)$-good if and only if $\rho = 1$.

The hyperplane $H$ is $(\rho,\ep)$-acceptable if and only if $\rho = 1$ or $\rho =\ep$ or $\rho$ is
the reflection representation or $n \leq 5$ or $n=6$ and $\rho$ is associated to one of the partition $(3,3)$, $(2,2,2)$, $(4,2)$.
\end{proposition}

\subsection{The rank 2 imprimitive groups}

For $d,e,r$ nonzero integers, we define the group $G(de,e,r)$ to be the group of  $r$-dimensional monomial matrices (with one nonzero element on each row and column)
whose nonzero entries are $de^{\ins{th}}$ root of unity such that their product is a $d^{\ins{th}}$ root of unity. These groups are called the imprimitive groups
of reflection. The integer $r$ is called the rank of $G(de,e,r)$.

The reflections of $G(de,e,2)$ are of the form
$$\begin{matrice} \xi &  \\ & 1 \end{matrice}\qquad\qquad \begin{matrice} 1 & \\ & \xi \end{matrice}\qquad\ins{and}\qquad
            \begin{matrice}  & \dz \\ \dz^{-1} &  \end{matrice}$$
\noindent where $\xi$ is a non trivial $d^{\ins{th}}$ root of unity and $\dz$ a $de^{\ins{th}}$ root of unity.
If $e$ is odd, there are two conjugacy classes of hyperplanes : one given by the hyperplanes of the diagonal reflections, the other given
by the nondiagonal reflections. If $e$ is even, there are three conjugacy classes of hyperplanes.
One given by the hyperplanes of the diagonal reflections. The hyperplanes of nondiagonal reflections split into two classes :
one associated to the hyperplane of
$$s= \begin{matrice}  & 1 \\ 1 &  \end{matrice}$$
\noindent the other one associated to the hyperplane of
$$s'= \begin{matrice}  & \dz \\ \dz^{-1} &  \end{matrice}$$
\noindent where $\dz$ is a $de^{\ins{th}}$ primitive root of unity.

We will describe the character of $G(de,e,2)$ following the method of small groups of Wigner and Mackey (see~\cite[paragraph 8.2]{serre}).
So we set $D$ the subgroup of diagonal matrix of $G(de,e,2)$. This is an abelian normal subgroup of $G(de,e,2)$ of index $2$.
To use the method of Wigner and Mackey, we have to describe the one-dimensional character of $D$.

\begin{lemma}[Linear character of $D$]\label{prop-caractere-lineaire-ddeer} For $d,e \in \NN^*$, the map
$$\fonction{\De}{\ZZ/de\ZZ \times \ZZ/d\ZZ}{\wh{D}=\Hom{gr}{D}{\CC^\times}}{(k,k')}{(\diag(\al,\be) \mapsto \al^{-k}(\al\be)^{-k'})}$$
\noindent is a group isomorphism.
\end{lemma}

Let us now describe the irreducible representation of $G(de,e,2)$. We have to distinguish with the evenness of $d$ and $e$.

\begin{proposition}[${d,e}$ odd]\label{prop-repres-de-impair} Let $d,e \in \NN^*$ be odd numbers.

For $k' \in \ZZ/d\ZZ$, we extend $\De(0,k')$ to $G(de,e,2)$ by $\De(0,k')(dx)= \De(0,k')(d)$ for every $d \in D$
and $x \in \langle 1,s \rangle$. We extend any irreducible representation $\rho$ of $\langle 1,s\rangle$ to $G(de,e,2)$ by $\rho(dx) = \rho(x)$
for every $d \in D$ and $x \in \langle 1,s \rangle$.
We then define
$$\be_{k',\rho}(dx) = \rho(dx)\De(0,k')(dx) = \rho(x)\De(0,k')(d)$$
for $d \in D$ and $x \in \langle 1,s\rangle$.
For $(k,k') \in \intnn{1}{de/2} \times \ZZ/d\ZZ$, we set $\be_{k,k'} = \ins{Ind}_D^{G(de,e,2)}\De(k,k')$.

The family $((\be_{k', 1},\be_{k', \ep})_{k' \in \ZZ/d\ZZ}, (\be_{k,k'})_{(k,k') \in \intnn{1}{de/2} \times \ZZ/d\ZZ})$
is a complete set for the irreducible representations of $G(de,e,2)$.
\end{proposition}

\begin{proposition}[${e}$ odd, ${d}$ even]\label{prop-repres-eimpair-dpair} Let $d,e \in \NN^*$ with $d=2d'$ even and $e$ odd.

For $k' \in \ZZ/d\ZZ$, we extend $\De(0,k')$ to $G(de,e,2)$ by $\De(0,k')(dx)= \De(0,k')(d)$ for every $d \in D$
and $x \in \langle 1,s \rangle$. We extend any irreducible representation $\rho$ of $\langle 1,s\rangle$ to $G(de,e,2)$ by $\rho(dx) = \rho(x)$
for every $d \in D$ and $x \in  \langle 1,s \rangle$. We then define
$$\be_{k',\rho}(dx) = \rho(dx)\De(0,k')(dx) = \rho(x)\De(0,k')(d)$$
for every $d \in D$ and $x \in \langle 1,s \rangle$.

We denote by $A$ the set $A=\paa{\intnn{1}{d'e-1} \times \ZZ/d\ZZ} \cup \paa{d'e} \times \intnn{0}{d'-1}$.
For $(k,k') \in A$, we set
$$\be_{k,k'} = \ins{Ind}_{D}^{G(de,e,2)}\De(k,k').$$

The family $((\be_{k', 1},\be_{k', \ep})_{k' \in \ZZ/d\ZZ}, (\be_{k,k'})_{(k,k') \in A})$
is a complete set of irreducible representations of $G(de,e,2)$.
\end{proposition}

\begin{proposition}[${e}$ even]\label{prop-repres-epair} Let $d,e \in \NN^*$ with $e=2e'$ even.

For $k' \in \ZZ/d\ZZ$ and $\de \in \paa{0,de'}$, we extend the character $\De(\de,k')$ to $G(de,e,2)$ by $\De(\de,k')(dx)= \De(\de,k')(d)$ for
$d \in D$ and $x \in \langle 1,s \rangle$.
We extend any irreducible representation $\rho$ of $\langle 1,s\rangle$ to $G(de,e,2)$ by $\rho(dx) = \rho(x)$ for every $d \in D$
and $x \in \langle 1, s \rangle$. We then define, for $d \in D$ and $x \in \langle 1,s \rangle$,
$$\be_{\de, k',\rho}(dx) = \rho(dx)\De(\de,k')(dx) = \rho(x)\De(\de,k')(d)\,.$$

\noindent For $(k,k') \in \intnn{1}{de'-1} \times \ZZ/d\ZZ$, we set
$$\be_{k,k'} = \ins{Ind}_{D}^{G(de,e,2)}\De(k,k')\,.$$

\noindent The family $((\be_{\de, k', 1},\be_{\de, k', \ep})_{(\de,k') \in \paa{0,de'} \times \ZZ/d\ZZ}, (\be_{k,k'})_{(k,k') \in \intnn{1}{de'-1} \times \ZZ/d\ZZ})$
is a complete set of irreducible representations of $G(de,e,2)$.
\end{proposition}

\begin{corollary}[Non diagonal hyperplanes] Let $d,e \in \NN^*$, $\rho$ be an irreducible representation of $G(de,e,2)$,
$\chi$ a linear character of $G(de,e,r)$ and $H$ the hyperplane of a non-diagonal reflection.

The hyperplane $H$ is $\rho$-good, $\rho$-excellent and $(\rho, \chi)$-acceptable for every $\rho$ and $\chi$.

If $e$ is odd and $\rho$ is a $2$-dimensional representation of $G(de,e,2)$,
the hyperplane $H$ is $(\rho,\chi)$-good if and only if $\chi = \be_{k',1}$ with $k' \in \ZZ/d\ZZ$.

If $e$ is odd, $\rho=\be_{k',\rho'}$ is a $1$-dimensional representation of $G(de,e,2)$ and $\chi = \be_{k'',\rho''}$ a linear character of $G(de,e,2)$,
then the hyperplane $H$ is $(\rho,\chi)$-good if and only if $\rho' \neq \ep$ or $\rho'' \neq \ep$.

Assume that $e=2e'$ is even, $H$ is an hyperplane associated to the conjugacy class of $s$ and $\rho$ a $2$-dimensional representation of $G(de,e,2)$
then $H$ is $(\rho,\chi)$-good if and only if $\chi = \be_{\de ,k',1}$ with $\de \in \paa{0,de'}$ and $k' \in \ZZ/d\ZZ$.

Assume that $e=2e'$ is even, $H$ is an hyperplane associated to the conjugacy class of~$s$, $\rho = \be_{u,k',\rho'}$ is a $1$-dimensional representation of $G(de,e,2)$
and $\chi = \be_{v ,k'',\rho''}$ a linear character of $G(de,e,2)$, then $H$ is $(\rho,\chi)$-good if and only if $\rho'' \neq \ep$ or $\rho' \neq \ep$.

Assume that $e=2e'$ is even, $H$ is an hyperplane associated to the conjugacy class of~$s'$ and $\rho$ a $2$-dimensional representation of $G(de,e,2)$
then $H$ is $(\rho,\chi)$-good if and only if $\chi = \be_{0 ,k',1}$ or $\chi = \be_{de',k',\ep}$ with $k' \in \ZZ/d\ZZ$.

Assume that $e=2e'$ is even, $H$ is an hyperplane associated to the conjugacy class of~$s'$, $\rho = \be_{u,k',\rho'}$ is a $1$-dimensional representation of $G(de,e,2)$
and $\chi = \be_{v ,k'',\rho''}$ a linear character of $G(de,e,2)$, then $H$ is $(\rho,\chi)$-good if and only if $(u,\rho') \notin \paa{(0,\ep), (de',{1})}$ or
$(v,\rho'') \notin \paa{(0,\ep), (de',{1})}$.
\end{corollary}

\begin{corollary}[Diagonal hyperplanes] Let $d,e \in \NN^*$ and $\be_{k,k'}$ be a $2$-dimensional irreducible representation of $G(de,e,2)$,
$\chi$ a linear character of $G(de,e,r)$ and $H$ an hyperplane associated to a diagonal reflection.

For $n \in \ZZ$, we denote by $\overline{n}$ the unique integer such that $0 \leq \overline{n} \leq d-1$ and $d \mid (n- \overline{n})$.

With these notations, we obtain $n_H(\be_{k,k'})= \overline{k+k'}+ k'$. For $e$ odd and $u \in \paa{1, \ep}$, we have $n_H(\be_{k',u}) = k'$
For $e=2e'$ even, $\de \in \paa{0,de'}$ and $u \in \paa{1, \ep}$, we have $n_H(\be_{\de,k', u}) = k'$.

If $k' = 0$ then $H$ is $\be_{k,k'}$-excellent and so $\be_{k,k'}$-good and $(\be_{k,k'},\chi)$-acceptable for every $\chi$.
Moreover $n_H(\be_{k,k'})=  \overline{k}$. Thus, for $\chi =\be_{\de,k'', u}$ or $\chi= \be_{k'',u}$ with $\de \in \paa{0,de'}$ and $u \in \paa{1, \ep}$,
the hyperplane $H$ is $(\be_{k,k'},\chi)$-good if and only if  $\overline{k}+k'' <d$.

Let us assume $k' \neq 0$. For $\chi =\be_{\de,k'', u}$ or $\chi= \be_{k'',u}$ with $\de \in \paa{0,de'}$ and $u \in \paa{1, \ep}$, the hyperplane $H$ is

\begin{ri}{(iii)}
\item $\be_{k,k'}$-excellent if and only if $\overline{k+k'}=0$;

\item $\be_{k,k'}$-good if and only if $\overline{k+k'} + k'< d$;

\item $(\be_{k,k'},\chi)$-good if and only if $\overline{k+k'} + k'+k'' < d$;

\item $(\be_{k,k'},\chi)$-acceptable if and only if $d-k'' > \overline{k+k'}$ or $d-k'' > k'$.
\end{ri}

If $e$ is odd, $\rho=\be_{k',\rho'}$ is a $1$-dimensional representation of $G(de,e,2)$ and $\chi = \be_{k'',\rho''}$ a linear character of $G(de,e,2)$,
then  $H$ is $(\rho,\chi)$-good if and only if $k'+k'' < d$.

Assume that $e=2e'$ is even, $\rho = \be_{u,k',\rho'}$ is a $1$-dimensional representation of $G(de,e,2)$
and $\chi = \be_{v ,k'',\rho''}$ a linear character of $G(de,e,2)$, then the hyperplane $H$ is $(\rho,\chi)$-good if and only if  $k'+k'' < d$.
\end{corollary}

\subsection{The wreath product}

Let us now study the imprimitive reflection group $G(d,1,r)$ which is also the wreath product $\ZZ/d\ZZ \wr \Sfr_r$.
Let us assume $r\geq 3$. The reflections of $G(d,1,r)$ are of the form
$$\diag(1,\ldots, 1, \xi, 1, \ldots, 1)$$
with $\xi$ a non trivial $d^{\ins{th}}$ of unity and of the form
$$\diag(1,\ldots,1,\dz, 1, \ldots, 1, \dz^{-1}, 1,\ldots, 1)\tau_{ij}$$
where $\tau_{ij}$ is the transposition matrix swapping $i$ and $j$ and $\dz$ is a
$d^{\ins{th}}$ root of unity. The associated hyperplanes split into two conjugacy class : the diagonal one and the non-diagonal one.

The irreducible character of $G(d,1,r)$ can be described by the method of Wigner and Mackey with the normal abelian subgroup of diagonal matrices of $G(d,1,r)$.
Thus the representation of $G(d,1,r)$ are given by the $d$-multipartitions of $r$ that is to say families of $d$ partitions so that the sum of the length of the
partitions are $r$. One can also describe the irreducible representations of $G(d,1,r)$ by giving a family of $d$ integers
$\underline{r} =(n_0,\ldots, n_{d-1})$ such that $n_0 + \cdots + n_{d-1} = r$ and $\rho$ a representation of $\Sfr_{n_0} \times \cdots \times \Sfr_{n_{d-1}}$.
We denote by $\be_{\underline{r},\rho}$ the corresponding representation of $G(d,1,r)$.

\begin{corollary}[Non diagonal hyperplanes]\label{cor-repres-gd1r-type-hyperplan}
Let $d \geq 2$, $r \geq 3$, $H$ be an non diagonal hyperplane of $G(d,1,r)$
(we denote by $G_H$ the subgroup of $G$ of reflections whose hyperplane is $H$), $\rho'=\be_{\underline{r},\rho}$ be
an irreducible representation of $G(d,1,r)$ and $\chi$ a linear character of $G(d,1,r)$.

\begin{ri}{(iii)}
\item The hyperplane $H$ is $\rho'$-excellent, $\rho'$-good, $(\rho',\chi)$-good (for $\chi(G_H)\!=\!1$),
$(\rho',\chi)$-acceptable (for $\chi(G_H)\!\neq\!1$) if and only if $\rho' = \be_{\underline{r},\rho}$ is of the form
\begin{rn}{a)}
\item[$a)$] $\underline{r} = (0, \ldots, 0, r, 0, \ldots, 0)$ and $\rho= 1$ or $\rho=\ep$ or $\rho$ is the standard representation
or  $\rho=(2,2)$ if $r=4$.

\item[$b)$] $\underline{r}= (0, \ldots, 0, 1, 0 \ldots, 0, r-1, 0, \ldots, 0)$ and $\rho=1$.

\item[$c)$] $\underline{r}= (0, \ldots, 0, r-1, 0 \ldots, 0, 1, 0, \ldots, 0)$ and $\rho=1$.
\end{rn}

\item $H$ is $(\rho',\chi)$-good (for $\chi(G_H) \neq 1$) if and only if $\underline{r} = (0, \ldots, 0, r, 0, \ldots, 0)$ and $\rho= 1$.

\item The hyperplane $H$ is $(\rho',\chi)$-acceptable (for $\chi(G_H) = 1$) if and only if $\rho'$ is one the following representation
\begin{rn}{a)}
\item[$a)$] $\underline{r} = (0, \ldots, 0, r, 0, \ldots, 0)$ and $\rho= 1$ or $\rho=\ep$ or $\rho$ is the standard representation
or $r \leq 5$ or $\rho \in \paa{(3,3), (2,2,2), (4,2)}$ if $r=6$.

\item[$b)$] $\underline{r}= (0, \ldots, 0, 1, 0 \ldots, 0, r-1, 0, \ldots, 0)$ and $\rho=1$ or $\rho=\ep$ if $r \in \paa{3,4}$
or $\rho$ is the standard representation if $r=3$.

\item[$c)$] $\underline{r}= (0, \ldots, 0, r-1, 0 \ldots, 0, 1, 0, \ldots, 0)$ and $\rho=1$ or $\rho=\ep$ if $r \in \paa{3,4}$
or $\rho$ is the standard representation if $r=3$.

\item[$d)$] $\underline{r}= (0, \ldots, 0, 2, 0 \ldots, 0, 3, 0, \ldots, 0)$ and $\rho=1$.

\item[$e)$] $\underline{r}= (0, \ldots, 0, 3, 0 \ldots, 0, 2, 0, \ldots, 0)$ and $\rho=1$.

\item[$f)$] $\underline{r}= (0, \ldots, 0, 2, 0 \ldots, 0, 2, 0, \ldots, 0)$ and $\rho=1$ or $\rho = \ep \otimes 1$ or $\rho = 1 \otimes \ep$.

\item[$g)$] $\underline{r}= (0, \ldots, 0, 1, 0 \ldots, 0, 1, 0, \ldots, 0, 1, 0, \ldots, 0)$ and $\rho = 1$.
\end{rn}
\end{ri}
\end{corollary}

\begin{corollary}[Diagonal hyperplanes] Let $d \geq 2$, $r \geq 3$, $H$ be a diagonal hyperplane of $G(d,1,r)$
and $\rho'=\be_{\underline{r},\rho}$ be an irreducible representation of $G(d,1,r)$ with $\underline{r}=(n_0,\ldots, n_{d-1})$
and $\chi$ a linear character of $G(d,1,r)$.

\begin{ri}{(iii)}
\item The hyperplane $H$ is $\rho'$-good if and only if $\ \ \dim\rho\,\pa{\Sum{j=0}{d-1}\frac{jn_j}{r}\,\frac{r!}{n_0! \cdots n_{d-1}!}}<d$
\item The hyperplane $H$ is $\rho'$-excellent if and only if $\rho'$ is one of the following representation
\begin{rn}{a)}
\item[$a)$] $\underline{r} = (r, 0, \ldots, 0)$;

\item[$b)$] $\underline{r} = (0, 0, \ldots, 0, r, 0, \ldots, 0)$ and $\rho=1$ or $\rho=\ep$;

\item[$c)$] $\underline{r} = (r-1, 0, \ldots, 0, 1, 0, \ldots, 0)$ and $\rho=1$ or $\rho=\ep$.
\end{rn}
\item The hyperplane $H$ is $(\rho',\chi)$-good if $\ \ \dim\rho\,\pa{\Sum{j=0}{d-1}\frac{jn_j}{r}\,\frac{r!}{n_0! \cdots n_{d-1}!}}<d-n_H(\chi)$.
\end{ri}
\end{corollary}

\subsection{The group ${G_{4}}$}

The group $G_{4}$ named after the Shephard and Todd classification~\cite{shephard-todd} is a rank $2$ reflection group.
There is only one class of hyperplanes which is of order $3$. The linear character of $G_{4}$ are given by the trivial one, the determinant and the
square of the determinant. As an abstract group, $G_{4}$ is nothing else but $SL(2,\FF_3)$.
The irreducible representations of $G_{4}$ are then given by the $3$ one-dimensional representation, the standard reflection representation named $V$,
and two others $2$-dimensional representation $V\det $ (whose character is real) and $V\det^2$ and one $3$-dimensional representation
(whose character is real).

\begin{corollary}[Hyperplanes of ${G_{4}}$]
Let $H$ be an hyperplane of $G_{4}$, $\rho$ be an irreducible representation of $G_4$.

Then $H$ is $\rho$-excellent, $\rho$-good, $(\rho,\det)$-acceptable and $(\rho,1)$-good
when $\dim \rho=1$ or $\rho = V$ or $\rho = V\det^2$.

The hyperplane $H$ is $(\rho,\det)$-good if and only if $\rho = 1$,

The hyperplane $H$ is $(\rho,\det^2)$-good if and only if $\rho = 1$, $\rho = \det^2$ or $\rho = V\det^2 $.

The hyperplane $H$ is $(\rho,1)$-acceptable and $(\rho,\det^2)$-acceptable for every $\rho$.
\end{corollary}

\subsection{The group ${G_{5}}$}

The group $G_{5}$ named after the Shephard and Todd classification~\cite{shephard-todd} is a rank $2$ reflection group.
There are two classes of hyperplanes which are both of order $3$. In fact,
$$G_5 = \{j^kG_{4},\ j = \exp(2i\pi/3), \quad k \in \{0,1,2 \}\} = G_4 \times \{ \id, j\id, j^2\id \}.$$
The irreducible representations of $G_{5}$ are then given by tensor product of representation of $G_4$ and of the three one-dimensional representation
of $\{ \id, j\id, j^2\id \}$ which are given by $(1,\det,\det^2)$.
One class of hyperplanes is in fact the class of hyperplanes of $G_4$. The other one is a new class.

\begin{corollary}[The hyperplanes of ${G_{5}}$ which are in ${G_4}$]
Let $H$ be an hyperplane of $G_{5}$ which is an hyperplane of $G_4$,
$\rho$ be an irreducible representation of $G_5$ and $\chi$ a linear character of $G_5$.
Write $\rho = \rho_1 \otimes \rho_2$ (resp. $\chi = \chi_1\otimes \chi_2$) where $\rho_1$ (resp. $\chi_1$) is an (resp. one-dimensional) irreducible representation of $G_4$
and $\rho_2$ (resp. $\chi_2$) an irreducible representation of $\{ \id, j\id, j^2\id \}$

Then $H$ is $\rho$-excellent if and only if $H$ is $\rho_1$-excellent for $G_4$.

The hyperplane $H$ is $\rho$-good if and only if $H$ is $\rho_1$-good for $G_4$.

The hyperplane $H$ is $(\rho,\chi)$-acceptable if and only if $H$ is $(\rho_1,\chi_1)$-acceptable for $G_4$.

The hyperplane $H$ is $(\rho,\chi)$-good if and only if $H$ is $(\rho_1,\chi_1)$-good for $G_4$.
\end{corollary}

\begin{corollary}[The hyperplanes of ${G_{5}}$ which are notin ${G_4}$]
Let $H$ be an hyperplane of $G_{5}$ which is not an hyperplane of $G_4$,
$\rho$ be an irreducible representation of $G_5$ and $\chi$ a linear character of $G_5$.

The hyperplane $H$ is $\rho$-excellent if and only if $H$ is $\rho$-good if and only if $\dim \rho = 1$, $\rho =V \otimes 1$,
$\rho =V \otimes \det^2$, $\rho=V\det \otimes \det$, $\rho=V\det \otimes \det^2$, $\rho=V\det^2 \otimes 1$, $\rho=V\det^2 \otimes \det$.

If $\chi=\det^{i} \otimes \det^{k}$ with $i+k = 0[3]$, then
$H$ is $(\rho,\chi)$-good if and only if $\dim \rho = 1$, $\rho =V \otimes 1$, $\rho =V \otimes \det^2$, $\rho=V\det \otimes \det$
$\rho=V\det \otimes \det^2$, $\rho=V\det^2 \otimes 1$, $\rho=V\det^2 \otimes \det$.

If $\chi=\det^{i} \otimes \det^{k}$ with $i+k = 1[3]$,
then $H$ is $(\rho,\chi)$-good if and only if $\rho=\det^{i'} \otimes \det^{k'}$ with $i'+k'=0[3]$ or $i'+k'=1[3]$ or $\rho =V \otimes 1$ or
$\rho=V\det \otimes \det^2$ or $\rho=V\det^2 \otimes \det$.

If $\chi=\det^{i} \otimes \det^{k}$ with $i+k = 2[3]$, then $H$ is $(\rho,\chi)$-good if and only if $\rho=\det^{i'} \otimes \det^{k'}$ with $i'+k'=0[3]$.

If $\chi=\det^{i} \otimes \det^{k}$ with $i+k = 0[3]$ or $i+k=1[3]$, the hyperplane $H$ is $(\rho,\chi)$-acceptable for every $\rho$.

If $\chi=\det^{i} \otimes \det^{k}$ with $i+k = 2[3]$,
the hyperplane $H$ is $(\rho,\chi)$-acceptable if and only if $\dim \rho = 1$, $\rho =V \otimes 1$, $\rho =V \otimes \det^2$, $\rho=V\det \otimes \det$
$\rho=V\det \otimes \det^2$, $\rho=V\det^2 \otimes 1$, $\rho=V\det^2 \otimes \det$.
\end{corollary}

\subsection{The group ${G_{24}}$}

The group $G_{24}$ named after the Shephard and Todd classification~\cite{shephard-todd} is a rank $3$ reflection group.
There is only one class of hyperplanes which is of order $2$. The linear character of $G_{24}$ are given by the determinant and the trivial one.
As an abstract group, $G_{24}$ is nothing else but the product of the simple groups $GL(3,\FF_2) \times \{ -1,1 \}$.
Let us denote by $1$ and $\ep$ the irreducible representations of $\{ -1,1 \}$ and
$1,3_1,3_2,6,7,8$ the irreducible representations of $GL(3,\FF_2)$ (determined by their dimension).
The irreducible representations of $G_{24}$ are then given by the tensor products of an irreducible representation of $GL(3,\FF_2)$
and $\{ -1,1 \}$.

\begin{corollary}[Hyperplanes of ${G_{24}}$]
Let $H$ be an hyperplane of $G_{24}$.

The hyperplane $H$ is $\rho$-excellent, $\rho$-good, $(\rho,1)$-good and $(\rho,\det)$-acceptable if and only if
$\rho = 1 \otimes 1$, $\rho = 1 \otimes \ep$, $\rho = 3_1 \otimes \ep$ and $\rho = 3_2 \otimes \ep$.

The hyperplane $H$ is $(\rho,\det)$-good if and only if $\rho= 1  \otimes 1$.

The hyperplane $H$ is $(\rho,1)$-acceptable if and only if $\rho = 1 \otimes 1$, $\rho = 3_1 \otimes 1$, $\rho = 3_2 \otimes 1$ and $\rho = 6 \otimes 1$ and
$\rho = 1 \otimes \ep$, $\rho = 3_1 \otimes \ep$ and $\rho = 3_2 \otimes \ep$ and $\rho = 7 \otimes \ep$.
\end{corollary}

\end{document}